\documentclass[11pt]{article}
\usepackage{amsmath}
\usepackage{amsfonts}
\usepackage{amssymb}
\usepackage{amsthm}
\usepackage{extarrows}
\usepackage{graphicx}
\usepackage{subfigure}

\usepackage{comment}

\usepackage{enumerate}
\newtheorem*{theorem*}{Lemma}
\newtheorem{theorem}{Theorem}
\newtheorem{lemma}{Lemma}
\newtheorem{corollary}[theorem]{Corollary}
\newtheorem{example}[theorem]{Example}
\newtheorem{lemma*}[theorem]{Lemma}
\newtheorem{definition}[theorem]{Definition}
\newtheorem{proposition}[theorem]{Proposition}

\begin{document}

\title{Orthonormal rational functions on a semi-infinite interval}

\author{Jianqiang Liu \footnote{School of Mathematics and Statistics, Ningxia University, Yinchuan, 750021, China. Email: liujq@amss.ac.cn} }

\maketitle
\begin{abstract}
  In this paper we propose a novel family of weighted orthonormal rational functions on a semi-infinite interval. We write a sequence of integer-coefficient polynomials in several forms and derive their corresponding differential equations. These equations do not form Sturm-Liouville problems. We overcome this disadvantage by multiplying some factors, resulting to a sequence of irrational functions. We deduce various generating functions of this sequence of irrational functions and find its associated Sturm-Liouville problems, which brings orthogonality. Then we study a Hilbert space of functions defined on a semi-infinite interval with its inner product induced by a weight function determined by the Sturm-Liouville problems mentioned above. We list two bases, one is the even subsequence of the irrational function sequence above and another one is the non-positive integer power functions.  We raise one example of Fourier series expansion and one example of interpolation as applications.
\end{abstract}

\noindent{\bf Keywords:}
Chebyshev polynomial $\cdot$ rational function  $\cdot$ orthonormal basis  $\cdot$ semi-infinite interval  $\cdot$ Hilbert space

\medskip

\noindent{\bf 2020 AMS Subject Classifications:}
 MSC 33D45,  41A20, 41A30,  46E30, 12D99

\newcommand{\C}[0]{\mathbb C}
\newcommand{\N}[0]{\mathbb N}
\newcommand{\R}[0]{\mathbb R}
\newcommand{\Z}[0]{\mathbb Z}
\section{Introduction}
Fibonacci sequence  is familiar to mathematical researchers (see \cite{frontczak2023balancing}, \cite{marques2022proof}, \cite{suvarnamani2015some}). Kilmer et. al. studied a kind of generalized Fibonacci sequence with a parameter,  namely \emph{$G$-polynomials}, through a three-term recurrence relation(see \cite{kilmer2023constructing}).
 
In Combinatorics (see \cite{brualdi2004introductory}), recurrence relation is an essential tool for generating sequences, for example, Fibonacci sequence and Chebyshev polynomial sequences(see \cite{szego1975orthogonal}). We are motivated to put forward the study $G$-polynomials from three-term recurrence relation, aiming at finding properties similar (or related) to that of Fibonacci sequence's and Chebyshev polynomial sequences'.  We start from writing $G$-polynomials in different forms and finding their corresponding differential equations. 
 
  The wide-speediness of Chebyshev Polynomials in various practical engineering circumstances relies a lot on the extensive results of its basic mathematical analysis properties. Orthogonality is their key character. Do $G$-polynomials have orthogonality?  We are interested in this question. Since the orthogonality of Chebyshev polynomials is embodied in a Sturm-Liouville problem(see \cite{asmar2016partial}), we investigate the differential equations that $G$-polynomials solve. The result has two sides. The up side is, there exist second order differential equations that $G$-polynomials solve. The down side is, in a standard form, the coefficient of first derivative of the unknown in each of these equations is a parameterized number, rather than a free one. This fact prevents $G$-polynomials from being an orthogonal sets according to Sturm-Liouville theory. 

Hence the family of $G$-polynomials is not the desired answer for the orthogonality. For the second best, we seek for a modification to take its place. This is accomplished by multiplying a factor to each $G$-polynomial, and the result is a sequence of irrational functions. Some properties of this sequence are derived, including generating functions and associated Sturm-Liouville problems.

Sturm-Liouville form brings orthogonality. By composing a Lebesgue measure and the weight function defined in Sturm-Liouville problems above, we establish a Hilbert space of functions on a semi-infinite interval. We find two of its bases, both are rational function sequences. One of them two is an orthonormal basis, it is the even subsequence of the irrational sequence mentioned above. Another one is not orthogonal, it is the function sequence of non-positive powers, which, after the Gram-Schmidt process, becomes the first basis. Some functions are raised as examples in this space. As applications of this Hilbert space and the orthonormal basis, one example of Fourier series expansion and one example of interpolation are given.

The rest of the article is arranged as follows. In Section \ref{sec: 2}, we list several forms of $G$-polynomials, some identities and related differential equations, and induce an irrational function sequence. In Section \ref{sec: 3},  we deduce three types of generating functions of this irrational sequence, and the associated Sturm-Liouville problems. As a highlight, in Section \ref{sec: 4}, we establish a function Hilbert space  and propose two bases, and give two examples as applications.

\section{Polynomial sequence $\{G_n\colon n\in \Z_+\}$}\label{sec: 2}
We first give some notations. Let $\Z$ be the set of all integers, $\Z_+$ be the set of all nonnegative integers, $\N$ be the set of all positive numbers. $\C$ be the set of all complex numbers, $\R$ be the set of all real numbers. Kilmer et.al. defined a sequence of polynomials, namely $G$-Polynomials, on complex domain in the following recursive way(see \cite{kilmer2023constructing}). We follow the notations there. Let both $G_0(z)$ and $G_1(z)$  be the constant $1$ and for each $n\ge 2$, let
\begin{equation}\label{eq1-01}
     \quad G_{n}(z) =G_{n-1}(z)-z G_{n-2}(z), \quad z\in \C.
\end{equation}
It can be written in closed form as 
\begin{equation}\label{eq1-02}
     G_{n}(z) = \sum_{k=0}^{\left
     \lfloor 
     \frac{n}{2}
     \right\rfloor}
     (-1)^k\binom{n-k}{k}z^k, \quad z\in \C, n\in \Z_+.
\end{equation}
In the following, we give several expressions of $G$-polynomials and the differential equations they correspond to. We also introduce an irrational function sequence from $G$-polynomials.
\subsection{Expressions}
 We first give a relation among $G$-polynomials and  Chebyshev polynomials of the second and the fourth kinds, $U_n(x)$ and $S_n(x)$. $U_n$ are defined originally on $[-1,1]$ by
\begin{equation}\label{eq61}
   U_n(x)=\frac{\sin(n+1)\arccos x}{\sin{\arccos x}}, \quad n\in\mathbb{Z}_+, \quad x \in (-1,1), 
\end{equation}
and then extended (via analytic extension) to be polynomials on $\mathbb{C}$. $S_n$ is derived from $U_n$ by 
\[
S_n(x) =U_n\left(\frac{x}{2}\right),\quad x\in\C, n\in\Z_+.
\]

The first ten of $G_n$'s and $S_n$'s are given below. Terms are arranged in reverse order. 
\begin{equation*}
    \begin{array}{ll}
    \begin{array}{l}
    G_0(z)  =  1,\cr
    G_1(z)  =  1,\cr
    G_2(z)  =  1- z,\cr
    G_3(z)  = 1 -2z,\cr
    G_4(z) =  1-3z+z^2,\cr
    G_5(z) = 1-4z+3z^2,\cr
    G_6(z) = 1-5z+6z^2-z^3,\cr
    G_7(z)= 1-6z+10z^2-4z^3,\cr
    G_8(z) = 1-7z+15z^2-10z^3+z^4,\cr
    G_9(z)= 1-8z+21z^2-20z^3+5z^4,
    \end{array}
         &  
    \begin{array}{l}
    S_0(z)  =  1,\cr
    S_1(z)  =  z,\cr
    S_2(z)  =  z^2-1 ,\cr
    S_3(z)  =  z^3-2z,\cr
    S_4(z) =  z^4-3z^2+1,\cr
    S_5(z) = z^5-4z^3+3z,\cr
    S_6(z) = z^6-5z^4+6z^2-1,\cr
    S_7(z) = z^7-6z^5+10z^3-4z,\cr
    S_8(z) = z^8-7z^6+15z^4-10z^2+1,\cr
    S_9(z) = z^9-8z^7+21z^5-20z^3+5z.
    \end{array}
    \end{array}
\end{equation*}
By comparison of the coefficients, we find for all $0\le k\le n$ and arbitrary $n\in \Z_+$, each coefficient of $z^k$ term in $G_n$ is the same as the coefficient of $z^{n-k}$ term in $S_n$. We prove this rule in the following theorem.
\begin{theorem}\label{thm7}
Suppose $z$ belongs to $\C^*$. Then for each $n\in\Z_+$, the following identity holds true:
\begin{equation} \label{eq01}
  G_n(z) = z^{\frac{n}{2}} S_n \left(\frac1{\sqrt{z}}\right) = z^{\frac{n}{2}} U_n \left(\frac1{2\sqrt{z}}\right).
\end{equation}
\end{theorem}
\begin{proof}
It suffices to prove the first equality. 
This clearly holds when $n=0$ and $1$.  By induction, for all $z\in \C^*$ and $n\ge 2$, we have
\begin{eqnarray*}
 G_n(z) &=& G_{n-1}(z) - z G_{n-2}(z) \nonumber\cr
 &=& z^{\frac{n-1}{2}} S_{n-1} \left(\frac1{\sqrt{z}}\right) - z z^{\frac{n-2}{2}} S_{n-2} \left(\frac1{\sqrt{z}}\right)\nonumber\cr
 &=& z^{\frac{n}{2}} \left[ \frac 1{\sqrt{z}} S_{n-1}\left(\frac1{\sqrt{z}}\right) - S_{n-2} \left(\frac1{\sqrt{z}}\right)\right] \nonumber\cr
 &=& z^{\frac{n}{2}} S_{n} \left(\frac1{\sqrt{z}}\right).
 \end{eqnarray*}
\end{proof}

For each $n \in \Z_+$, consider the two-variable function $F_n(z,x)$ defined by
\begin{equation*}
  F_n(z,x)=z^{\frac n2}U_n\left (\frac{x}{2\sqrt{z}}\right), \quad x \in [-1,1],  z \in \C.
\end{equation*}
We see that $G_n(z)=F_n(z,1)$ for all $z\in \mathbb{C}$, i.e., $G_n$ is a special case of $F_n(z,x)$  with the second variable $x$ fixed at $1$. Kiepiela and  Klimek \cite{kiepiela2005extension} studied the case of the same function family $F_n(z,x)$ when the first variable $z$ is fixed.

If $z=0$, we observe that
\begin{equation}\label{eq03}
    G_n(0)=1\quad \textrm{for all } n\in\Z_+.
\end{equation}
At $z=\frac{1}{4}$, the analyticity of $G_n$ enables us to establish the following combinatorial formula.
\begin{corollary}\label{cor01}
For each $n\in \mathbb{Z}_+$,
 \begin{equation}\label{eq10}
    G_n\left(\frac{1}{4}\right)=\sum_{k=0}^{\lfloor \frac {n}{2} \rfloor}(-1)^k\binom{n-k}{k}\frac{1}{4^k}=\frac{n+1}{2^n}.
  \end{equation}
\end{corollary}
\begin{proof}
 If $z=\frac 14$, we may write $G_0$ and $G_1$ in the following form:
  \[
  G_0\left(\frac{1}{4}\right)=1=\frac{0+1}{2^0},\quad G_1\left(\frac{1}{4}\right)=1=\frac{1+1}{2^1}.
  \]
By induction, let's suppose 
 \[
 G_k\left(\frac{1}{4}\right)=\frac{k+1}{2^k} \quad \textrm{ for } 0\le k <n.
 \]
Then we have
 \begin{equation*}
G_n\left(\frac{1}{4}\right)=\frac{n}{2^{n-1}}-\frac{1}{4} \frac{n-1}{2^{n-2}}=\frac{n+1}{2^n} , \quad n\in \mathbb{Z}_+.
\end{equation*}
The proof completes since each $G_n$ is analytic on $\C$. 
\end{proof}

Theorem \ref{thm7} allows us to rewrite $G_n$ in terms of trigonometric functions.
\begin{corollary}\label{cor05}
   Suppose $z$ belongs to $\C^* \setminus\{\frac{1}{4}\}$. For each $n\in \mathbb{Z}_+$, there holds
\begin{equation}\label{eq02}
G_n(z)=
\frac{2z^{\frac{n+1}{2}}}{\sqrt{4z-1}}\sin(n+1)\arccos \frac{1}{2\sqrt{z}}.
\end{equation}
\end{corollary}
\begin{proof}
It can be shown by substitution $z=\frac{1}{4\cos^2 \theta}$ in Theorem 1 in Reference \cite{kilmer2023constructing}. We also give a proof using Lemma 1 in this paper. When $z\ne 0$ and $z\ne \frac{1}{4}$,
  \begin{eqnarray*}
    G_n(z)&=&z^{\frac{n}{2}}U_n\left(\frac{1}{2\sqrt{z}}\right) \cr
    &=&\frac{z^{\frac{n}{2}}\sin(n+1)\arccos \left(\frac{1}{2\sqrt{z}}\right)}{\sin\arccos \left(\frac{1}{2\sqrt{z}}\right)}\cr
    &=&\frac{2z^{\frac{n+1}{2}}}{\sqrt{4z-1}}\sin(n+1)\arccos \left(\frac{1}{2\sqrt{z}}\right).
  \end{eqnarray*}
\end{proof}

Hence $G_n$ is the analytical extension of the right hand side of \eqref{eq02} in Corollary \ref{cor05} to entire complex plane after adding  $G_n(0)$ and $G_n\left(\frac{1}{4}\right)$ as defined in \eqref{eq03} and \eqref{eq10} respectively. Now we give another expression of $G_n$. In Combinatorics, a three term recurrence relation such as \eqref{eq1-01} corresponds to a \emph{characteristic equation} (see \cite{brualdi2004introductory}). By definition, the characteristic equation of  \eqref{eq1-01} is 
\begin{equation}\label{eq07}
    s^2-s+z=0,
\end{equation}
which is a quadratic equation with parameter $z$. It's two solutions are respectively
\begin{equation}\label{eq75}
    s_1(z) =\frac{1-\sqrt{1-4z}}{2} \quad \textrm{and}\quad s_2(z) =\frac{1+\sqrt{1-4z}}{2}.
\end{equation}
In terms of $s_1(z)$ and $s_2(z)$, we can rewrite $G_n$ as the following form.
\begin{corollary}\label{cor02}
For each $n\in \Z_+$ and $z\in \C\setminus\{\frac{1}{4}\}$, 
\begin{equation} \label{eq14}
    G_n(z)= \frac{s_2^{n+1}(z)-s_1^{n+1}(z)}{s_2(z)-s_1(z)}.
\end{equation}
\end{corollary}
\begin{proof}
Suppose $z \ne \frac{1}{4}$. Clearly \eqref{eq14} holds true for $n=0$ and $n=1$. Now let $n\ge 2$. Suppose 
$$G_k(z) =\frac{s_2^{k+1}(z)-s_1^{k+1}(z)}{s_2(z)-s_1(z)} \quad \textrm{for}\quad  0\le k<n.$$ 
Then by \eqref{eq1-01},
\begin{eqnarray}
G_n(z) &=&\frac{s_2^{n}(z)-s_1^{n}(z)}{s_2(z)-s_1(z)} - z\frac{s_2^{n-1}(z)-s_1^{n-1}(z)}{s_2(z)-s_1(z)} \nonumber\\
&=&\frac{s_2^{n-1}(z)(s_2(z)-z)- s_1^{n-1}(z)(s_1(z)-z)}{s_2(z)-s_1(z)}. \label{eq16}
\end{eqnarray}
Since both $s_1(z)$ and $s_2(z)$ are roots of equation \eqref{eq07}, we obtain 
\[
s_j(z) -z =s_j^2(z), \quad j =1,2.
\]
This enables us to rewrite $G_n$ from \eqref{eq16} to
\[
G_n(z) = \frac{s_2^{n-1}(z)s_2^2(z)- s_1^{n-1}(z)s_1^2(z)}{s_2(z)-s_1(z)} = \frac{s_2^{n+1}(z)-s_1^{n+1}(z)}{s_2(z)-s_1(z)}.
\]
Then we ends our proof by induction.
\end{proof}

\medskip

Note that if $z =-1$, we have 
\[
s_1(-1) = \frac{1-\sqrt{5}}{2}, \quad s_2(-1) =\frac{1+\sqrt{5}}{2}.
\]
Hence we can write $G_n(-1)$ in the form of \eqref{eq14} as
\begin{equation}\label{eq74}
    G_n(-1) =\frac{1}{\sqrt{5}}\left[\left(\frac{1+\sqrt{5}}{2}\right)^{n+1}-\left(\frac{1-\sqrt{5}}{2}\right)^{n+1}\right].
\end{equation}

Recall that Fibonacci sequence is defined by $F_0=F_1=1, F_{n}= F_{n-1} +F_{n-2}$ for $n\ge 2$, possessing a formula of general term of equation \eqref{eq74}. Hence we have 
\begin{corollary}\label{lm14}
For all $n\in \Z_+$, $F_n =G_n(-1)$.
\end{corollary}

\subsection{Differential Equations}

It is known that the second kind of Chebyshev polynomials $U_n$ solves equation
  \begin{equation}\label{eq05}
    (1-x^2)\frac{d^2v}{dx^2}-3x\frac{dv}{dx}+n(n+2)v=0,
  \end{equation}
 This motivates us to find the differential equation related to each polynomial $G_n$. It is shown in the following theorem.
\begin{theorem}\label{thm01}
  For each $n\in\mathbb{Z}_+$, $v=G_n(z)$ is the solution to the following equation:
  \begin{equation} \label{eq31}
  (4z^2-z)v'' +[n-(4n-6)z]v'  +n(n-1)v=0.
\end{equation}
\end{theorem}
\begin{proof}
    If $x\ne 0$, by substitution $x=\frac{1}{2\sqrt{z}}$ in equation \eqref{eq01}, we have $z=\frac{1}{4x^2}$ and
 \begin{equation*}
   U_n(x)=(2x)^n G_n\Big(\frac{1}{4x^2}\Big).
 \end{equation*}
Then we derive from \eqref{eq05} that 
\begin{equation}\label{eq82}
  \Big(\frac{1}{4 x^{4}}-\frac{1}{4 x^{2}}\Big)G''_n\Big(\frac{1}{4x^2}\Big)
  +\left[
  n-\Big(\frac{2n-3}{2}\Big)\frac{1}{x^{2}}
  \right]G'_n\Big(\frac{1}{4x^2}\Big)
  +n(n-1)G_n\Big(\frac{1}{4x^2}\Big)=0.
\end{equation}
Now we write the left hand side of \eqref{eq82} in terms of $z$ again. This gives
\begin{equation*}
  (4z^2-z)G''_n(z) +[n-(4n-6)z]G'_n(z)  +n(n-1)G_n(z)=0. 
\end{equation*}
\end{proof}

It can be seen from equation \eqref{eq31} that the coefficient of $v'$ term is relative to $n$, which prevent $G_n$ from being a solution to a Sturm-Liouville problem. For the second best, we try to find a certain function  related $G_n$, namely $\Phi_n$, coming to stage as a solution to a related Sturm-Liouville problem. In specific operation, we seek for some factor function $\phi_n(z)$, such that $\Phi_n$, defined by $\Phi_n(z) :=\phi_n(z)G_n(z)$, is a solution of second order ordinary equation in which the coefficient of $\Phi'_n$ terms, denoted by $\alpha(n,z)$, is 
unrelated to $n$ provided the the coefficient of $\Phi''_n$ term is $1$. To achieve this, we substitute $v = \Phi_n(z)\phi_n^{-1}(z)$ into \eqref{eq31}, and obtain that
\begin{equation}\label{eq83}
    (4z^2-z)\left(\frac{\Phi_n(z)}{\phi_n(z)}\right)'' +[n-(4n-6)z]\left(\frac{\Phi_n(z)}{\phi_n(z)}\right)'  +n(n-1)\left(\frac{\Phi_n(z)}{\phi_n(z)}\right)=0.
\end{equation}
Equation \eqref{eq83} may then be transformed to the following form
\begin{eqnarray}
  &&\frac{4z^2-z}{\phi_n^{3}}
  [\phi_n^2\Phi''_n-2 \phi_n\phi'_n\Phi'_n +(2\phi'^2_n-\phi_n\phi'_n)\Phi_n] \nonumber\\
  &&\quad \quad \quad   +\frac{n-(4n-6)z}{\phi_n^{2}}(\phi_n\Phi'_n-\phi'_n\Phi_n)+\frac{n(n-1)}{\phi_n}\Phi_n =0. \label{eq84}
\end{eqnarray}
By rearranging terms in equation \eqref{eq84}, we obtain
\begin{eqnarray}
  &&\frac{4z^2-z}{\phi_n}\Phi''_n + 
\left[\frac{n-(4n-6)z}{\phi_n}-\frac{2(4z^2-z)}{\phi_n^2}\phi'_n\right]\Phi'_n
+\frac{1}{\phi_n^3}\{
(4z^2-z)\nonumber\\ 
&&\quad\quad\quad (2\phi'^2_n-\phi_n\phi'_n)+[n-(4n-6)z]\phi_n+n(n-1)\phi_n^2
\}\Phi_n =0. \label{eq33}
\end{eqnarray}
When $z\ne 0$ and $z\ne \frac{1}{4}$, let us divide equation \eqref{eq33} by $\frac{4z^2-z}{\phi_n}$, we get
\begin{eqnarray}
  &&\Phi''_n + 
  \left[\frac{n-(4n-6)z}{4z^2-z} -\frac{2\phi'_n}{\phi_n}\right]\Phi'_n \label{eq49}\\
&+&\frac{(4z^2-z)(2\phi'^2_n-\phi_n\phi'_n)+[n-(4n-6)z]\phi_n+n(n-1)\phi_n^2}{\phi_n^2(4z^2-z)}\Phi_n =0. \nonumber
\end{eqnarray}
We deduce from \eqref{eq49} that 
\begin{equation}\label{eq50}
    \alpha(n,z)=\frac{n-(4n-6)z}{4z^2-z}-2\frac{\phi'_n}{\phi_n}
    =\frac{6}{4z-1}-\frac{n}{z}-2\frac{\phi'_n}{\phi_n}.
\end{equation}
To eliminate $n$ in $\alpha(n,z)$ in \eqref{eq50}, a simple choice is to let
\begin{equation}\label{eq41}
    -\frac{n}{z}-2\frac{\phi'_n}{\phi_n} =0.
\end{equation}
Equation \eqref{eq41} is a first order ordinary differential equation. Its general solution is 
\begin{equation*}
    \phi_n(z) =\frac{C}{z^{\frac{n}{2}}}, \quad C \textrm{    is a constant.}
\end{equation*}
That is to say, to make the coefficient of each $\Phi_n$'s unrelated with $n$,
we can choose $\Phi_n$  as 
\begin{equation*}
    \Phi_n(z) =\frac{CG_n(z)}{z^{\frac{n}{2}}}.
\end{equation*}
By letting $C=1$, we obtain a particular solution and denote it by $g_n$. See the definition below. 
\begin{definition}\label{def1}
For each $n\in \Z_+$ and $z\in \C^*$, define function $g_n$ as follows:
\begin{equation}\label{eq43}
    g_n(z)=\sum_{k=0}^{\left
     \lfloor 
     \frac{n}{2}
     \right\rfloor}
     (-1)^k\binom{n-k}{k}z^{k-\frac{n}{2}}.
\end{equation}
\end{definition}
From Definition \ref{def1} we see that each $g_n$ is an irrational function, and each $g_{2n}$ is a rational function. In the next section we will focus on  $\{g_n\colon n\in Z_+\}$. We will call $\{g_{2n}\colon  n \in Z_+\}$ the set of \emph{$g$-rational functions}. We will give more about its properties in the Section \ref{sec: 4}.

\section{Irrational function sequence $\{g_n\colon n\in Z_+\}$}\label{sec: 3}
In this section, we characterize irrational function sequence $\{g_n\colon n\in Z_+\}$, from its basic properties, generating functions and corresponding differential equations. 
\subsection{Generating functions}
In this part we give three kinds of generating functions. We begin with the following lemma, which reveals the relation among $G_n,U_n$ and $g_n$, and  gives two types of recurrence relations.
\begin{lemma}\label{lm02}
  For each $n\in\Z_+$, let $g_n$ be defined as in \eqref{eq43}, then each $g_n(z)$ is analytic for all $z\in \C^*$, and the following identities hold true:
  \begin{eqnarray}
    g_n(z) &=&\frac{G_n(z)}{z^{\frac{n}{2}}},\label{eq66}\\
    g_n(z) &=& U_n\left(\frac{1}{2\sqrt{z}}\right),\label{eq67}\\
    g_{n+2}(z) &=& \frac{1}{\sqrt{z}}g_{n+1}(z)-g_{n}(z),\label{eq68}\\
    g_{n+4}(z) &=& \left(\frac{1}{z}-2\right)g_{n+2}(z)-g_{n}(z). \label{eq69}
  \end{eqnarray}
\end{lemma}
\begin{proof}
Analyticity follows from $g_n$'s definition. Identity \eqref{eq66} follows from equation \eqref{eq1-02}, Identity \eqref{eq67} follows from equation \eqref{eq01}. Identity \eqref{eq68} follows from  \eqref{eq66} and equation \eqref{eq1-01}. Now we prove Identity \eqref{eq69}. For each $n\in\Z_+$ and  $z\in \C^*$, making use of Identity \eqref{eq68}, we can deduce that
\begin{eqnarray*}
  g_{n+4}(z) &=& \frac{1}{\sqrt{z}}g_{n+3}(z)-g_{n+2}(z) \cr
  &=& \frac{1}{\sqrt{z}}\left( \frac{1}{\sqrt{z}}g_{n+2}(z)-g_{n+1}(z)\right)-g_{n+2}(z) \cr
  &=& \frac{1}{z}g_{n+2}(z)-\frac{1}{\sqrt{z}}g_{n+1}(z)-g_{n+2}(z) \cr
  &=& \frac{1}{z}g_{n+2}(z)-\left(g_{n+2}(z)+g_{n}(z)\right)-g_{n+2}(z) \cr
  &=& \left(\frac{1}{z}-2\right)g_{n+2}(z)-g_{n}(z).
\end{eqnarray*}
\end{proof}

In the following lemma, we evaluate $g_n$ at $z=\frac{1}{4}$ and $z=\infty$.
 \begin{lemma}\label{lm07}
  For each $n\in\Z_+$, let $g_n$ be defined as in \eqref{eq43}, then the following two identities hold true:
  \begin{eqnarray}
    g_n\left(\frac{1}{4}\right) &=&n+1,\label{eq70}\\
    g_n(\infty) &=& \left\{
  \begin{array}{ll}
     (-1)^{\frac{n}{2}},  & \textrm{if $n$ even,} \cr
     0,  & \textrm{if $n$ odd.}
    \end{array}
       \right.\label{eq71}
  \end{eqnarray}
\end{lemma}
\begin{proof}
Both identities can be derived from the property of $U_n$ and Identity \eqref{eq67}. Here we present a proof using definition of $g_n$ directly.
\begin{enumerate}[(i).]
      \item By equation \eqref{eq10} in Corollary \ref{cor01}, if $z= \frac{1}{4}$, we have 
\[
g_n\left(\frac{1}{4}\right) = 2^n G_n\left(\frac{1}{4}\right) = 2^n \cdot \frac{n+1}{2^n}= n+1.
\]
      \item When $z=\infty$, by equation \eqref{eq1-02},
\[
g_n(\infty)=\lim_{z\to \infty}\sum_{k=0}^{\left
     \lfloor 
     \frac{n}{2}
     \right\rfloor}
     (-1)^k\binom{n-k}{k}z^{k-\frac{n}{2}} =\sum_{k=0}^{\left
     \lfloor 
     \frac{n}{2}
     \right\rfloor}
     (-1)^k\binom{n-k}{k}\lim_{z\to \infty}z^{k-\frac{n}{2}}.
\]
Recall that $\lim_{z\to \infty}z^j =0$ if $j<0$. If $n$ is odd, then
\[
\frac{n}{2} >k\quad \textrm{when  }0\le k\le \left\lfloor 
     \frac{n}{2}
     \right\rfloor
\]
and $g_n(\infty)=0$. On the other hand, if $n$ is even, then we have
\[
g_n(\infty)=(-1)^\frac{n}{2}\binom{\frac{n}{2}}{\frac{n}{2}}\lim_{z\to \infty}z^{0}=(-1)^\frac{n}{2}.
\]
\end{enumerate}

\end{proof}

Now we define a semi-infinite interval as follows:
$$\mathcal{Z} = \left[\frac{1}{4}, +\infty\right),$$
Boyd established an orthogonal rational basis on a semi-infinite interval $(0, +\infty)$ in \cite{BOYD198763}. Now we give three kinds of generating functions of $\{g_n,n\in \Z_+\}$. 

\begin{proposition}\label{prop1}
For $z\in \mathcal{Z}$ and $|t|< 1$, the following identity hold true:
\begin{equation*}
\sum_{n=0}^\infty g_n(z) t^n = \frac{1}{1-z^{-\frac{1}{2}}t+ t^2}.
\end{equation*}
\end{proposition}
\begin{proof}
\noindent
We start with the well-known identity for Chebyshev polynomials of the second kind.
\[
\sum_{n=0}^\infty U_n(x) t^n= \frac{1}{1-2x t+t^2},  \quad t \in (-1, 1), \quad x \in [-1, 1].
\]
For $z\in \mathcal{Z}$, let $x=\frac{1}{2\sqrt{z}}$, 
then  $ x \in (0, 1]$. Hence for $z\in \mathcal{Z}$, $|t|< 1$,  we have
 \begin{equation*}
\sum_{n=0}^\infty g_n(z) t^n =\sum_{n=0}^\infty U_n\left(\frac{1}{2\sqrt{z}}\right) t^n 
= \frac 1{1-2\left(\frac{1}{2\sqrt{z}}\right)t + t^2} =\frac{1}{1-z^{-\frac{1}{2}}t+ t^2}. \label{eq13}
\end{equation*}
\end{proof}

\begin{proposition}\label{prop2}
For $z\in\mathcal{Z}$, $t\in \mathbb{R}$, the following identity holds true:
    \begin{equation*}
        \sum_{n=0}^\infty g_n(z) \frac{t^n}{n!} = \left\{ \begin{array}{ll}
 e^{\frac{t}{2\sqrt{z}}}\left(
  \cos \frac{t\sqrt{4z-1}}{2\sqrt{z}} +
  \frac{1}{\sqrt{4z-1}}\sin \frac{t\sqrt{4z-1}}{2\sqrt{z}}
  \right), & \textrm{if } z>\frac{1}{4}, \vspace{2mm}\cr
\displaystyle(t +1) e^{2t}, &\textrm{if } z= \frac14.
\end{array} 
\right.
    \end{equation*}
\end{proposition}
\begin{proof}
Let $S(z,t)$ be defined as 
\begin{equation}\label{eq76}
    S(z,t) = \sum_{n=0}^\infty g_n(z) \frac{t^n}{n!}, \quad z>\frac{1}{4}, t\in \R.
\end{equation}
Then 
\begin{equation}\label{eq77}
\frac{\partial S}{\partial t} = \sum_{n=1}^\infty g_n(z) \frac{t^{n-1}}{(n-1)!} = \sum_{n=0}^\infty g_{n+1}(z) \frac{t^n}{n!},
\end{equation}
and 
\begin{equation}\label{eq78}
\frac{\partial^2 S}{\partial t^2} =  \sum_{n=0}^\infty g_{n+2}(z) \frac{t^n}{n!}.
\end{equation}
Applying \eqref{eq76}, \eqref{eq77} and Identity \eqref{eq68}, we can rewrite $g_{n+2}$ in \eqref{eq78} to obtain
\begin{equation} \label{eq79}
    \frac{\partial^2 S}{\partial t^2} -\frac{1}{\sqrt{z}} \frac{\partial S}{\partial t} + S =0.
\end{equation}
Now consider the following algebraic equation with a parameter $z$:
\[
x^2 -\frac{x}{\sqrt{z}}  + 1 =0.
\]
Its solutions are respectively
\[
x_{1}=\frac{s_{1}(z)}{\sqrt{z}},\quad x_{2}=\frac{s_{2}(z)}{\sqrt{z}},
\]
where $s_1(z)$, $s_2(z)$ are defined as in equation \eqref{eq75}. Hence the general solution to equation \eqref{eq79} is 
  \[
  S(z,t)= C_1(z)e^{\frac{s_{1}(z)t}{\sqrt{z}}} + C_2(z)e^{\frac{s_{2}(z)t}{\sqrt{z}}},
  \]
where $C_1(z)$ and $C_2(z)$ are undetermined functions of $z$. Since 
\[
S(z,0) =1,\quad \frac{\partial S}{\partial t}(z,0) = \frac{1}{\sqrt{z}}, 
\]
we have the following linear system:
\begin{equation}\label{eq80}
\left\{
    \begin{array}{rcrcc}
       C_1(z)&+&C_2(z)  & \equiv & 1,  \\
        s_1(z)C_1(z)&+&s_2(z)C_2(z) & \equiv&1 . 
    \end{array}
    \right.
\end{equation}
The solutions of \eqref{eq80} are
\[
C_1(z) = \frac{is_1(z)}{\sqrt{4z-1}} ,\quad 
C_2(z) = -\frac{is_2(z)}{\sqrt{4z-1}}.
\]
Hence we have
\begin{eqnarray*}
  S(z,t) &=& \frac{i}{\sqrt{4z-1}}\left[s_1(z)e^{\frac{s_1(z)t}{\sqrt{z}}} - s_2(z)e^{\frac{s_2(z)t}{\sqrt{z}}}\right]\\
  &=& \frac{i e^{\frac{t}{2\sqrt{z}}}}{2\sqrt{4z-1}}
  \left[(1-i\sqrt{4z-1})
  \left(
  \cos \frac{t\sqrt{4z-1}}{2\sqrt{z}} -i
  \sin \frac{t\sqrt{4z-1}}{2\sqrt{z}}
  \right) 
  \right.
  \\
  &&
  \quad\quad\left.
  -
  (1+i\sqrt{4z-1})
  \left(
  \cos \frac{t\sqrt{4z-1}}{2\sqrt{z}} +i
  \sin \frac{t\sqrt{4z-1}}{2\sqrt{z}}
  \right)
  \right]\\
  &=&e^{\frac{t}{2\sqrt{z}}}\left(
  \cos \frac{t\sqrt{4z-1}}{2\sqrt{z}} +
  \frac{1}{\sqrt{4z-1}}\sin \frac{t\sqrt{4z-1}}{2\sqrt{z}}
  \right).
\end{eqnarray*}
If $z=\frac14$, \eqref{eq70} imples
\[
\sum_{n=0}^{\infty}\frac{g_n\left(\frac14\right)t^n}{n!} 
= \sum_{n=0}^{\infty}\frac{(n+1)t^n}{n!} =(t+1)e^t.
\]
\end{proof}

\begin{proposition}\label{prop3}
For $z\in \mathcal{Z}$ and $|t|< 1$, the following identity holds true:
\[
\sum_{n=1}^\infty g_n(z) \frac{t^n}n =\left\{ 
\displaystyle\begin{array}{ll} 
\displaystyle\frac{1}{\sqrt{4z-1}}\arctan \frac{2\sqrt{z}t-1}{\sqrt{4z-1}} - \frac{1}{2}\ln|1- z^{-\frac{1}{2}}t +t^2|, & \textrm{if } z> \frac14,\vspace{2mm}\cr
 \displaystyle \frac{1}{1-t}-\ln(1-t), & \textrm{if }  z=\frac14.
\end{array}
\right.
\]
\end{proposition}

\begin{proof}
Since
\[
\frac{\partial}{\partial t}\left(
\sum_{n=1}^{\infty}g_n(z)\frac{t^n}{n}\right)= \frac1t\left(\sum_{n=0}^\infty g_n(z) t^n  - g_0(z)\right),
\]
we deduce that 
\begin{eqnarray}
  \sum_{n=1}^{\infty} g_n(z)\frac{t^n}{n}
  &=&\int_{0}^t\frac{1}{\tau}\left(
  \frac{1}{1-z^{-\frac{1}{2}}\tau +\tau^2}  -1
  \right)\, d\tau\nonumber\\
  &=& \int_{0}^t \frac{z^{-\frac{1}{2}}-\tau}{1-z^{-\frac{1}{2}}\tau +\tau^2}\, d\tau\nonumber\\
  &=&\frac{1}{2}\int_{0}^t\left(
  \frac{z^{-\frac{1}{2}}-2\tau}{1-z^{-\frac{1}{2}}\tau +\tau^2}
  +\frac{z^{-\frac{1}{2}}}{1-z^{-\frac{1}{2}}\tau +\tau^2}
  \right)\, d\tau. \label{eq81}
\end{eqnarray}
If $z>\frac{1}{4}$, \eqref{eq81} gives
\[
\sum_{n=1}^{\infty} g_n(z)\frac{t^n}{n} =\frac{1}{\sqrt{4z-1}}\arctan \frac{2\sqrt{z}t-1}{\sqrt{4z-1}} - \frac{1}{2}\ln|1- z^{-\frac{1}{2}}t +t^2|.
\]
If $z=\frac{1}{4}$,  \eqref{eq81} implies
\[
\sum_{n=1}^{\infty} g_n\left(\frac{1}{4}\right)\frac{t^n}{n} =\frac{1}{1-t}-\ln(1-t).
\]
\end{proof}

\medskip
\subsection{Sturm-Liouville problems}
For each $z\in \mathcal{Z}$, let 
\begin{equation}\label{eq26}
    \theta= \arctan \sqrt{4z-1}, \quad \theta \in \left[0,\frac{\pi}{2}\right).
\end{equation}
Then this identity is a bijection between $\mathcal{Z}$ and $\left[0,\frac{\pi}{2}\right)$.  In the following lemma we connect $g_n$ with $\theta$ and the well-known $Dirichlet$ kernels(see \cite{rudin1976principles}). We recall that for each $n\in\Z_+$, the $n$-th Dirichlet kernel, denoted by $D_n$, is defined as follows
\[
D_n(\theta) = \frac{\sin \left(n+\frac{1}{2}\right)\theta }{\sin\frac{\theta}{2}}, \quad  \theta \in [-\pi, \pi].
\] 

\begin{lemma}\label{lm12}
  Let $\theta \in \left(0,\frac{\pi}{2}\right)$ be defined as in \eqref{eq26}. For each $n\in\Z_+$, let $g_n$ be defined as in \eqref{eq43}. Then for all $z>\frac{1}{4}$,  the following two identities hold true:
  \begin{eqnarray}
    g_n(z) &=& \displaystyle \frac{\sin (n+1)\theta}{\sin \theta},\label{eq72}\\
    g_{2n}(z) &=&D_n(2 \theta).\nonumber
  \end{eqnarray}
\end{lemma}
\begin{proof}
\begin{enumerate}[(i)]
    \item Since $\tan \theta =\sqrt{4z-1}$, we have $\cos \theta = \frac{1}{2\sqrt{z}}$. By Corollary \ref{cor05} and  Lemma \ref{lm02},
\[
g_n(z) =\frac{\sin(n+1)\arccos\frac{1}{2\sqrt{z}}}{\sin\arccos \frac{1}{2\sqrt{z}}} =\frac{\sin(n+1)\arctan\sqrt{4z-1}}{\sin\arctan\sqrt{4z-1}} =\frac{\sin(n+1)\theta}{\sin\theta}.
\]
\item By equation \eqref{eq72} we have 
\[
g_{2n}(z)  = 
\frac{\sin(2n+1)\theta}{\sin\theta}=\frac{\sin\left[\left(n+\frac{1}{2}\right)\cdot2\theta\right]}{\sin\frac{2\theta}{2}} =D_n(2\theta).
\]
\end{enumerate}
\end{proof}

Now we calculate the derivative of $g_n$ at $\frac{1}{4}$ for each $n\in\Z_+$.
\begin{lemma}\label{lm08}
  For each $n\in\Z_+$, let $g_n$ be defined as in \eqref{eq43}, then there holds 
  \begin{equation}\label{eq29}
      g'_n\left(\frac{1}{4}\right) =-\frac{2}{3}n(n+1)(n+2).
  \end{equation}
\end{lemma}
\begin{proof}
Since for all $z\in \mathcal{Z}$,
     $$g_0(z)\equiv1,\quad g_1(z)=\frac{1}{\sqrt{z}},$$ 
     we have  
     $$g'_0(z)\equiv 0,\quad g'_1(z)=-\frac{1}{2z^\frac{3}{2}}.$$
     Let $z=\frac{1}{4}$,  we obtain
     $$g'_0\left(\frac{1}{4}\right)=0,\quad g'_1\left(\frac{1}{4}\right)=-4.$$
     Now suppose \eqref{eq29} holds for all $m<n$. From identities \eqref{eq66} and \eqref{eq68} in Lemma \ref{lm02}, we have for all $n>1$ and $z\in \mathcal{Z}$,
     \begin{equation}\label{eq30}
         g'_n(z) =  \frac{1}{\sqrt{z}}g'_{n-1}(z)-g'_{n-2}(z)-\frac{1}{2z^\frac{3}{2}}g_{n-1}(z),
     \end{equation}
     hence 
       \begin{eqnarray*}
         g'_n\left(\frac{1}{4}\right) & = &2g'_{n-1}\left(\frac{1}{4}\right)-g'_{n-2}\left(\frac{1}{4}\right)-4g_{n-1}\left(\frac{1}{4}\right) \cr
         &=& -\frac{4}{3}n(n-1)(n+1) + \frac{2}{3}n(n-1)(n-2)-4n \cr
         &=& -\frac{2}{3}n(n-1)(n+1) \cr
         &&+\left[-\frac{2}{3}n(n-1)(n+1)+ \frac{2}{3}n(n-1)(n-2)-4n \right]\cr
         &=& -\frac{2}{3}n(n-1)(n+1) -2n(n+1)\cr
         &=& -\frac{2}{3}n(n+1)(n+2).
       \end{eqnarray*}
       By induction, \eqref{eq29} holds for all $n\in\Z_+$.
\end{proof}

\begin{lemma}\label{lm09}
  For each $n\in\Z_+$, let $g_n$ be defined as in \eqref{eq43}. Then  $g'_n(+\infty)=0$.
\end{lemma}
\begin{proof}
It is clear that 
$$g'_0(+\infty)=0$$
and
\[
g'_1(+\infty) := \lim_{z\to +\infty}-\frac{1}{2z^\frac{3}{2}} =0.
\]
 Now fix $n>1$ and suppose 
 \[
 g'_{n-1}(+\infty)=0 \quad \textrm{and}\quad g'_{n-2}(+\infty)=0.
 \]
 By Lemma \ref{lm07}, each $g_n$ is bounded at $+\infty$. Hence
 \[
 g'_n(+\infty) = \lim_{z\to +\infty} \frac{1}{\sqrt{z}}g'_{n-1}(z)-\lim_{z\to +\infty}g'_{n-2}(z)-\lim_{z\to +\infty}\frac{1}{2z^\frac{3}{2}}g_{n-1}(z) = 0.
 \]
 Hence $g'_n(+\infty)=0$ for each $n\in\Z_+$.
\end{proof}

The construction of $g_n$ enables us to establish the following theorem about a singular Sturm-Liouville problem on $\mathcal{Z}$(see \cite{asmar2016partial}).

\begin{theorem}\label{thm03}
Suppose that for each $n \in \Z_+, z \in \mathcal{Z}$, $g_n(z)$ is defined as in equation \eqref{eq43}. Set 
\begin{equation*}
    p(z) =(4z-1)^{\frac{3}{2}}, \quad q(z)=0, \quad  r(z) = \frac{\sqrt{4z-1}}{4z^2}, \quad \lambda =n(n+2).
\end{equation*}
Then $g_n$ is a solution to the following singular Sturm-Liouville problem
\begin{eqnarray}\label{eq45}
  &&[p(z)v']'+[q(z)+\lambda r(z)]v = 0, \quad z\in \left(\frac{1}{4}, +\infty\right),\cr
 &&\left\{\begin{array}{lcl}
       v\left(\frac{1}{4}\right) &=& n+1,  \cr
       v(+\infty)&=&
       \left\{
  \begin{array}{ll}
      (-1)^{\frac{n}{2}},  & \textrm{if $n$ even,}\cr
     0,  & \textrm{if $n$ odd.}
    \end{array}
       \right.
  \end{array}
  \right.
\end{eqnarray}
\end{theorem}
\begin{proof}
For $z\in \mathcal{Z}$, we substitute $v=g_n(z)$ into the left hand side of equation \eqref{eq45}, and deduce that
\begin{eqnarray*}
  && [p(z)g_n']'+[q(z)+\lambda r(z)]g_n \cr
  &=&\left[(4z-1)^{\frac{3}{2}} \left(\frac{G_n}{z^{\frac{n}{2}}}\right)'\right]' + \frac{n(n+2)\sqrt{4z-1}G_n}{4z^{\frac{n}{2}+2}} \cr
  &=&\left[\frac{2z(4z-1)^{\frac{3}{2}}G_n'-n(4z-1)^{\frac{3}{2}}G_n}{2z^{\frac{n}{2}+1}}\right]' + \frac{n(n+2)\sqrt{4z-1}G_n}{4z^{\frac{n}{2}+2}} \cr
  &=&\frac{4z^2(4z-1)^{\frac{3}{2}}G''_n + 4z[n-(4n-6)z]\sqrt{4z-1}G'_n
  }{4z^{\frac{n}{2}+2}} \cr
  &&+\frac{[4n(n-1)z-n(n+2)]\sqrt{4z-1}G_n
  }{4z^{\frac{n}{2}+2}} 
  + \frac{n(n+2)\sqrt{4z-1}G_n}{4z^{\frac{n}{2}+2}} \cr
  \medskip\cr
  &=&\frac{\sqrt{4z-1}
  }{z^{\frac{n}{2}+1}} \{(4z^2-z)G''_n + [n-(4n-6)z]G'_n+n(n-1)\sqrt{4z-1}G_n\}\cr
  &=&0.
\end{eqnarray*}

Note that the last step applies equation \eqref{eq31} in Theorem \ref{thm01}. The conditions of determining solution follows immediately from Lemma \ref{lm07}.
\end{proof}

By Lemmas \ref{lm08} and \ref{lm09}, the conditions of determining solution in \eqref{eq45} can be replaced by the derivatives at $z=\frac{1}{4}$ and $z=+\infty$, which is shown in the following lemma
\begin{corollary}
For each $n\in \Z_+$, $g_n$ solves the following equation:
\begin{eqnarray*}
  &&(16z^3-4z^2)v''+24z^2v'+n(n+2)v=0, \quad z \in \left(\frac{1}{4}, +\infty\right),\cr
 &&\left\{\begin{array}{lcl}
       v'\left(\frac{1}{4}\right) &=& -\frac{2}{3}n(n+1)(n+2),  \cr
       v'(+\infty)&=& 0.
  \end{array}
  \right.
\end{eqnarray*}
\end{corollary}

Now both the coefficients of $v''$ and $v'$ are unrelated to parameter $n$.
With Sturm-Liouville theory, we can discuss the orthogonality of function sequence $\{g_{2n}\colon n\in \Z_+\}$  on $\mathcal{Z}$ in the coming section.

\section{$L^2(\mathcal{Z},\mathcal{M},\mu)$ space}\label{sec: 4}

In this section we discuss the orthogonality related to $\{g_n\colon n\in \Z_+\}$ in a Hilbert space. In the first part we give some basic definitions of this space, and in the second part we give the orthonormal basis for it.
\subsection{Basic settings}
In this part we sequentially give several definitions including a measure, spaces, norms, an inner product and generalized Fourier series. Now let $\mathcal{M}$ be the $\sigma$-field consisting of all Lebesgue measurable subsets of $\mathcal{Z}$. Then $(\mathcal{Z},\mathcal{M})$ is a measurable space. Let $\mu$ be a set function defined by
	\begin{equation*}
	\mu(E) := \int_{E}\frac{\sqrt{4z-1}}{2\pi z^2}\, dz, \quad E \in \mathcal{M}.
	\end{equation*}
 For more about measure theory, please refer to textbooks of real analysis, such as  \cite{royden2010real}.  We have the following lemma.
\begin{lemma}
	$\mu$ is a probability measure on $(\mathcal{Z},\mathcal{M})$.
\end{lemma}
\begin{proof}
	Clearly $\mu(E)$ belongs to  $[0, +\infty]$ for all $E \in \mathcal{M}$. For empty set $\Phi$,
	\begin{equation*}
	\int_{\Phi}\frac{\sqrt{4z-1}}{2\pi z^2}\, dz=0.
	\end{equation*}
	Now let $\{E_k\colon k\ge 1\} \subset \mathcal{M}$. We have 
	\begin{eqnarray*}
		\mu(\bigcup_{k\ge 1}E_k)=\int_{\bigcup_{k\ge 1}E_k}\frac{\sqrt{4z-1}}{2\pi z^2}\, dz
		=\sum_{k=1}^{\infty}\int_{E_k}\frac{\sqrt{4z-1}}{2\pi z^2}\, dz=\sum_{k=1}^{\infty}\mu(E_k),
	\end{eqnarray*}
	hence $\mu$ is a measure on measure space $(\mathcal{Z},\mathcal{M})$. \eqref{eq26} implies
	\begin{equation}\label{eq51}
	    d\mu =  \frac{\sqrt{4z-1}}{2\pi z^2}\, dz 
	=  \frac{\tan \theta}{2\pi}16\cos^4\theta  \cdot  \frac{1}{2}\sec^2\theta \tan \theta \, d\theta 
	=  \frac{4}{\pi}\sin^2\theta \, d\theta.
	\end{equation}
	 Since $\theta\in (0,\frac{\pi}{2}]$ when $z\in \mathcal{Z}$, we have 
	\begin{equation*}
	\mu(\mathcal{Z}) =  \int_{\mathcal{Z}}\, d\mu 
	=  \int_{0}^{\frac{\pi}{2}}\frac{4}{\pi}\sin^2\theta \, d\theta =1, 
	\end{equation*}
	hence $\mu$  is a probability measure on $(\mathcal{Z}, \mathcal{M})$.
\end{proof}
\medspace

A set $E$ is said to be $\mu$-measurable if $E \in \mathcal{M}$. The readers can check that $\mu$ is absolutely continuous with respect to Lebesgue measure on $\mathcal{Z}$ and is a complete measure. A function $f$ on $\mathcal{Z}$ is said to be $\mu$-measurable if for each $c\in \mathbb{R}$, $\{x\colon f(x)>c\}$ is $\mu$-measurable. For more details, please refer to Royden and Fitzpatrick, \cite{royden2010real}.

For all $p\in [1,+\infty]$, we regard two functions $f$ and $g$ 
as one if they are equal on $\mathcal{Z}\setminus E$ with $\mu(E)=0$.
 For each $\mu$-measurable function $f$, define the $p$-norm of $f$, denoted by $\|f\|_p$, as follows:
\begin{equation*}
  \|f\|_p=\left\{
  \begin{array}{ll}
     \left(\int_{\mathcal{Z}}|f|^p\, d\mu\right)^{\frac{1}{p}},  &  \textrm{  if } p\in [1,+\infty),\\
      \mathrm{ess} \, \sup_{z\in\mathcal{Z}} |f(z)|, & \textrm{  if } p=+\infty, 
  \end{array}
  \right.
  \end{equation*}
where $\mathrm{ess} \, \sup |f|$ is the essential supremum of  $f$. Let $L^p(\mathcal{Z},\mathcal{M},\mu)$ be the collection of real-valued $\mu$-measurable functions on $\mathcal{Z}$ for which $\int_{\mathcal{Z}}|f|^p\, d\mu<\infty$.  Readers who are interested in more properties, such as H\"older inequality and Minkowski inequality, can refer to the textbooks of functional analysis(see \cite{conway2019course},\cite{stein2011functional}).
  
If $p=2$, the inner product of $f, g\in L^2(\mathcal{Z},\mathcal{M},\mu)$ is defined by
\begin{equation*}
(f,g)=\int_{\mathcal{Z}}fg\, d\mu,
\end{equation*}
then $L^2(\mathcal{Z},\mathcal{M},\mu)$ forms a Hilbert space. For each $f\in L^2(\mathcal{Z},\mathcal{M},\mu)$, the norm $\|f\|_2$ can be written as
\begin{equation*}
    \|f\|_2 = (f,f)^{\frac{1}{2}}.
\end{equation*}

Let $\mathcal{I}$ be an index set and $f_j\in L^p(\mathcal{Z},\mathcal{M},\mu)$ for each $j\in \mathcal{I}$. A collection $\mathcal{F}=\{f_j\colon j\in \mathcal{I}\}$ is said to be a (Schauder) \emph{basis} for $L^p(\mathcal{Z},\mathcal{M},\mu)$, if for each $f\in L^p(\mathcal{Z},\mathcal{M},\mu)$, and $\epsilon>0$,
there exists an integer $n$ and real numbers $c_k$ and functions $f_{j_k}\in L^p(\mathcal{Z},\mathcal{M},\mu), j_k\in \mathcal{I}$ for which $\|f-\sum_{k=0}^nc_kf_{j_k}\|_p<\epsilon$. If $p=2$, $\mathcal{F}$ is said to be \emph{an orthonormal set}, if for each $k,l\in \mathcal{I}$, there holds
\begin{equation*}
(f_k,f_l)=\delta_{k,l},
\end{equation*}
where $\delta_{kl}$ is the Kronecker delta defined by
$$\displaystyle  \delta_{jk}=\left\{
\begin{array}{lc}
1,& \textrm{if } k=l; \cr
0,& \textrm{otherwise.}
\end{array}
\right.$$

$\mathcal{F}$ is said to be an \emph{orthonormal basis} for $L^2(\mathcal{Z},\mathcal{M},\mu)$ if it is a maximal orthonormal set in $L^2(\mathcal{Z},\mathcal{M},\mu)$. For each $f\in L^2(\mathcal{Z},\mathcal{M},\mu)$ and  $f_j\in \mathcal{F}$, $(f,f_j)$ is called a (generalized) Fourier coefficient, $\sum_{j\in \mathcal{I}}(f,f_j)f_j$ is called the (generalized) Fourier series of $f$ with respect to $\mathcal{F}$ if there exists a function sequence $\{f_n\colon n\in \Z_+\}\subset \mathcal{F}$ and a unique real sequence $\{c_n\colon n\in \Z_+\}$ for which $\lim_{n\to \infty}\|f-\sum_{k=0}^nc_kf_k\|_2=0$. Then we may write $f=\sum_{j\in \mathcal{I}}(f,f_j)f_j$, and $\sum_{j\in \mathcal{I}}(f,f_j)f_j$ is called the Fourier expansion of $f$ with respect to basis $\mathcal{F}$. We have the following example.
\begin{example}
\begin{enumerate}[(i)]
\item Let $z_0>\frac{1}{4}$. Then  $(z-z_0)^\gamma \in L^2(\mathcal{Z},\mathcal{M},\mu)$ if and only if $\gamma \in (-\frac{1}{2},\frac{1}{4})$\vspace{2mm}.
    \item $(z-\frac{1}{4})^\gamma\in L^2(\mathcal{Z},\mathcal{M},\mu)$ if and only if $\gamma \in (-\frac{3}{4},\frac{1}{4})$.
\end{enumerate}
\end{example}
\begin{proof}
We prove Case (i). The proof of Case (ii) is similar. Suppose $\frac{1}{4}<z_0<z_1$. Then
\begin{equation*}
    \int_{\mathcal{Z}} (z-z_0)^{2\gamma}\, d\mu = \int_{\frac{1}{4}}^{+\infty} (z-z_0)^{2\gamma} \frac{\sqrt{4z-1}}{z^2}\, dz =I_1+I_2,
\end{equation*}
where
\[
I_1 =\int_{\frac{1}{4}}^{z_1} (z-z_0)^{2\gamma} \frac{\sqrt{4z-1}}{z^2}\, dz, \quad I_2=\int_{z_1}^{+\infty} (z-z_0)^{2\gamma} \frac{\sqrt{4z-1}}{z^2}\, dz.
\]
Improper integral $I_1$ converges if and only if  $2\gamma> -1$, i.e. $\gamma> -\frac{1}{2}$,  and $I_2$ converges if and only if $2\gamma +\frac{1}{2}+1<2$, i.e. $\gamma< \frac{1}{4}$. The range of $\gamma$ is then obtained by taking intersection.
\end{proof}

\subsection{An orthonormal basis}

 The following theorem is our main theorem in this section. 
\begin{theorem}\label{thm6}
$\{g_{2n}\colon  n \in Z_+\}$ forms an orthonormal basis for $L^2(\mathcal{Z},\mathcal{M},\mu)$.
\end{theorem}
To prove this theorem, we need several lemmas.
\begin{lemma}\label{lm06}
	$\{g_{2n}\colon  n \in Z_+\}$ forms an orthonormal set in $L^2(\mathcal{Z},\mathcal{M},\mu)$.
\end{lemma}
\begin{proof}
	For $m, n\in\mathbb{Z}_+$, we have
	\begin{eqnarray*}
	    (g_{2m}, g_{2n}) &=& \int_{\frac{1}{4}}^{+\infty}g_{2n}(z)g_{2m}(z)\frac{\sqrt{4z-1}}{2\pi z^2}\, dz \\
	    &=& \frac{4}{\pi}\int_{+\infty}^{\frac{1}{4}}g_{2n}(z)g_{2m}(z)\sqrt{1-\frac{1}{4z}} \cdot \left(-\frac{dz}{4z^{\frac{3}{2}}}\right).
	\end{eqnarray*}
	Let $z=\frac{1}{4x^2}, x \in (0, 1]$. Using Lemma \ref{lm02}, we obtain
	\begin{equation}\label{eq85}
	    (g_{2m}, g_{2n})  = \frac{4}{\pi}\int_{0}^1U_{2n}(x)U_{2m}(x)\sqrt{1-x^2}\, dx.
	\end{equation}
	Since $U_{2n}$ is even for each $n$, equation \eqref{eq85} implies 
	\begin{equation*}
	    (g_{2m}, g_{2n})  = \frac{2}{\pi}\int_{-1}^1U_{2n}(x)U_{2m}(x)\sqrt{1-x^2}\, dx = \delta_{2m,2n} =\delta_{m,n}.
	\end{equation*}
\end{proof}
This proof also shows, for each pair of $m, n\in \Z_+$, an inner product of $g_{2m}(z)$ and $g_{2n}(z)$ in $L^2(\mathcal{Z},\mathcal{M},\mu)$ can be transformed directly to be an inner product of $U_{2m}(x)$ and $U_{2n}(x)$ with respect to the weight function $\sqrt{1-x^2}$ through substitution $z=\frac{1}{4x^2}$ for $x \in [-1, 1]$.

Now we define
\[\mathcal{H} =  \overline{\mathrm{span}}\{g_{2n}\colon n\in \Z_+\}.\]
Then $\mathcal{H}$ is a closed subspace of $L^2(\mathcal{Z},\mathcal{M},\mu)$,  and $\{g_{2n}\colon n\in \Z_+\}$ is an orthonormal basis for $\mathcal{H}$ by a theorem about orthogonal basis(see \cite{conway2019course}). In the following we will prove that $\mathcal{H}=L^2(\mathcal{Z},\mathcal{M},\mu)$. Now set
\begin{equation*}
  C_1(\mathcal{Z}) =\{f\colon f \textrm{  is continuous on } \mathcal{Z} \textrm{  and } \lim_{z\to +\infty}f(z) \textrm{  exists}\}.
\end{equation*}
Then we have the following lemma.
\begin{lemma}\label{lm03}
  $C_1(\mathcal{Z}) \in \mathcal H$.
\end{lemma}
\begin{proof}
 Suppose $\lim_{z\to +\infty} f(z) =  a$. Let $z= \frac{1}{4x^2}$, then $z\to +\infty$ whenever $x\to  0^+$. This fact enables us to extend $f(\frac{1}{4x^2})$ to an even continuous function on $[-1, 1]$, namely $F(x)$, written as  
  \begin{equation*}
    F(x)=\left\{
      \begin{array}{ll}
        f(\frac{1}{4x^2}),& \hbox{if $x\in [-1,0)\bigcup(0, 1)]$, and }\cr
        a,& \hbox{if $x=0$}.
      \end{array}
      \right.
  \end{equation*}
  For any $n\in\mathbb{Z}_+$ and $c_k\in \R$, $0\le k\le n$, by applying the symmetry of even order $U_{2k}$, we deduce that
  \begin{eqnarray}
    \Big\|f-\sum_{k=0}^n c_kg_{2k}\Big\|_2^2 
    &=&\int_{\frac{1}{4}}^{+\infty} \big(f(z)-\sum_{k=0}^n c_kg_{2k}(z)\big)^2\frac{\sqrt{4z-1}}{2\pi z^2}\, dz\nonumber\\
    &=&\frac{4}{\pi}\int_{0}^1\sqrt{1-x^2}\Big(f(\frac{1}{4x^2})-\sum_{k=0}^n c_k U_{2k}(x)\Big)^2\, dx\nonumber\\
     &=&\frac{2}{\pi}\int_{-1}^1\sqrt{1-x^2}\Big(F(x)-\sum_{k=0}^n c_k U_{2k}(x)\Big)^2\, dx \nonumber\\
    &=& \Big\|F-\sum_{k=0}^n c_k U_{2k}\Big\|^2,\label{eq34}
  \end{eqnarray}
  where $\|\cdot\|$ denote the Euclidean norm with respect to the weight function $\sqrt{1-x^2}$ defined on $[-1, 1]$. By a corollary of Weierstrass's theorem(a version with respect to Chebyshev polynomials and $L^2$ norm on the interval $[-1,1]$, see Corollary 3.2A in \cite{mason2002chebyshev}), for given $\epsilon>0$, there exists $m\in\mathbb{Z}_+$ and $\{a_k\colon 0\le k \le m\}\subset \mathbb{R}$ such that 
  \begin{equation}\label{eq36}
    \Big\|F-\sum_{k=0}^ma_k U_{k}\Big\| < \frac{\epsilon}{2}.
  \end{equation}
  Since $F$ is even, by the substitution $t=-x$, we have
  \begin{eqnarray}
    \frac{\epsilon}{2}&>&\Big\|F-\sum_{k=0}^ma_k U_{k}\Big\| \nonumber\\
     &=&\Big(
      \frac{2}{\pi}\int_{-1}^1\sqrt{1-x^2}\Big(F(x)-\sum_{k=0}^ma_k U_{k}(x)\Big)^2\, dx  
    \Big)^{\frac{1}{2}}\nonumber\\
    &=&\Big(
      \frac{2}{\pi}\int_{-1}^1\sqrt{1-t^2}\Big(F(t)-\sum_{k=0}^m(-1)^k a_k U_{k}(t)\Big)^2\, dt  
    \Big)^{\frac{1}{2}}\nonumber\\
    &=&\Big\|F-\sum_{k=0}^m(-1)^k a_k U_{k}\Big\|.\label{eq39}
  \end{eqnarray}
Now we choose $n=\lfloor \frac{m}{2}\rfloor$ and $c_k=a_{2k}, 0\le k\le n$. Then similar to \eqref{eq39} we deduce that 
\begin{eqnarray}
  \Big\|f-\sum_{k=0}^{n}c_{k}g_{2k}\Big\|_2 &=&\Big\|F-\sum_{k=0}^{n}c_{k}U_{2k}\Big\| \nonumber\\
  &=& \Big\|\frac{1}{2}\Big(F-\sum_{k=0}^{m}a_{k}U_{k}\Big)+\frac{1}{2}\Big(F-\sum_{k=0}^{m}(-1)^k a_k U_{k}\Big)\Big\|.\label{eq17}
\end{eqnarray}
By applying Minkowski inequality we deduce from \eqref{eq36}-\eqref{eq17} that 
\begin{equation*}
  \Big\|f-\sum_{k=0}^{n}c_{k}g_{2k}\Big\|_2 \le \frac{\sqrt{2}}{2}\Big\|F-\sum_{k=0}^{m}a_{k}U_{k}\Big\|+\frac{\sqrt{2}}{2}\Big\|\Big(F-\sum_{k=0}^{m}(-1)^k a_k U_{k}\Big)\Big\|
    < \epsilon.
  \end{equation*}
 By the arbitrariness of $\epsilon$ and the completion of $\mathcal{H}$, we have $f\in \mathcal{H}$.
\end{proof}
We have the following examples.

\begin{example}
For each $n\in \Z_+$, each $g_{n}$ is continuous on $\mathcal{Z}$. By Lemma \ref{lm07},
\[
g_n(+\infty) =\left\{
\begin{array}{ll}
(-1)^\frac{n}{2},     & \textrm{ if $n$ even,}\\
0,     & \textrm{ if $n$ odd.}
\end{array}
\right.
\]
Hence $g_n \in  C_1(\mathcal{Z})\subset\mathcal{H}$  for $n\in \Z_+$.
 \end{example}
 
\begin{example}
Let $\beta , m \in \R$ and $f_{\beta,m}(z)=e^{\beta z^m}, z\in \mathcal{Z}$. Then 
\[
\lim_{z\to +\infty}  e^{\beta z^m}=\left\{
\begin{array}{ll}
0,     & \textrm{ if $m>0,\beta<0$, }  \\
e^{\beta},     & \textrm{ if $m=0$, }  \\
1,     & \textrm{ if $\beta =0$ or $m<0$,} 
\end{array}
\right.
\]
and $\lim_{z\to +\infty}  e^{\beta z^m}= +\infty$ if $\beta>0$ and $m>0$.
So we conclude that  $f_{\beta,m}\in \mathcal{H}$ iff $\beta m\le 0$(we follow the convention that $z^0 = 1$ at $z=0$. Otherwise, the condition "$\beta m\le 0$" should be replaced by "$\beta m \le 0$ and at least one of $\beta$ and $m$ is nonzero").
\end{example}
 
 \begin{example}\label{example1}
 Let $\deg(P)$ denote the degree of polynomial $P$. Let $f =\frac{P}{Q}$, where $P$ and $Q$ are relatively prime polynomials with $\deg(P)\le \deg(Q)$, and $Q$ has no zeros in $\mathcal{Z}$. Then
 $f$ is continuous on $\mathcal{Z}$, and 
 \[
\lim_{z\to +\infty}  f(z)=\left\{
\begin{array}{ll}
\displaystyle\frac{\alpha_P}{\alpha_Q},     & \textrm{ if $\deg(P) = \deg(Q)$,}  \\
0,     & \textrm{ if $\deg(P) < \deg(Q)$,}   
\end{array}
\right.
\]
where $\alpha_P$ and $\alpha_Q$ are the leading coefficients of $P$ and $Q$ respectively. Hence  $f \in \mathcal{H}$.
 \end{example}
 
 We continue searching for functions in $\mathcal{H}$.
 
 \begin{lemma}\label{lm11}
  Let $E$ be one of the following three cases:
 \begin{enumerate}[(i)]
     \item $E =[\frac{1}{4}, a)$, $a \in \mathcal Z$;
     \item $E$ is an open interval on $\mathcal{Z}$;
     \item $E\in \mathcal{M}$.
 \end{enumerate}
  Let $\mathbf{1}_E$ be the indicator function of $E$ for $E\in \mathcal{M}$, i.e.,
\begin{equation*}
    \mathbf{1}_E(z) = \left\{
\begin{array}{lc}
1, & \textrm{if } z \in E; \cr
0, & \textrm{otherwise.}
\end{array}
\right.
\end{equation*}
Then $\mathbf{1}_E \in \mathcal{H}$.
 \end{lemma}
 \begin{proof}
 \begin{enumerate}[(i)]
     \item Without loss of generality, we suppose $a>2$. For each $n\ge 1$ we define 
\begin{equation*}
    f_{n,a}(z) = \left\{
\begin{array}{ll}
1,& \displaystyle \textrm{if } z \in [\frac{1}{4},a-\frac{1}{n}], \cr
\displaystyle\frac{1}{2}(1+na-nz),& \displaystyle\textrm{if } z \in (a-\frac{1}{n}, a+\frac{1}{n}), \cr
0,& \textrm{otherwise.}
\end{array}
\right.
\end{equation*}
That is, each $f_{n,a}$ is identical with $\mathbf{1}_E$ in both $[\frac{1}{4},a-\frac{1}{n}]$ and $[a+\frac{1}{n}, +\infty)$, meanwhile draws a straight line segment connecting $(a-\frac{1}{n}, 1)$ and $(a+\frac{1}{n}, 0)$ for  $z\in (a-\frac{1}{n},a+\frac{1}{n})$. Then $f_{n,a}$ is continuous on $\mathcal{Z}$ and $\lim_{z\to +\infty}f_{n,a} (z)= 0$. By Lemma \ref{lm03}, we have $f_{n,a}\in \mathcal{H}$. Furthermore,  
\begin{eqnarray*}
    \|\mathbf{1}_E -f_{n,a}\|^2_2&=&\int_{a-\frac{1}{n}}^{a+\frac{1}{n}} (1-f_{n,a}(z))^2\frac{\sqrt{4z-1}}{2\pi z^2}\, dz \\
    &\le& \frac{2}{n} \max_{z\in (a-\frac{1}{n},a+\frac{1}{n})}\frac{\sqrt{4z-1}}{2\pi z^2} \le \frac{2}{\pi n} \to 0,\quad n\to \infty.
\end{eqnarray*}
Hence $\mathbf{1}_E \in \mathcal{H}$ for $E =[\frac{1}{4}, a)$.
\item Suppose $E= (b,c), c>b\ge \frac{1}{4}$. We have $\mathbf{1}_{\{z=b\}} \in \mathcal{H}$ since $\mathbf{1}_{\{z=b\}} = 0$ $\mu$-a.e. on $\mathcal{Z}$. Now applying the result of Case (i), we get
\begin{equation*}
    \mathbf{1}_{E} = \mathbf{1}_{[\frac{1}{4},c)} -\mathbf{1}_{[\frac{1}{4},b)}-\mathbf{1}_{\{z=b\}} \in \mathcal{H}.
\end{equation*}
Hence $\mathbf{1}_E \in \mathcal{H}$ for each open interval $E\subset \mathcal{Z}$.
\item Now given $\epsilon>0$ and $E\in \mathcal{M}$, there exist $z_1>\frac{1}{4}$ such that 
\begin{equation*}
    \mu([z_1,+\infty)) <\frac{\epsilon^2}{2}.
\end{equation*}
Let $E_0 =E\bigcap [\frac{1}{4}, z_1]$. Using a theorem of measure(Theorem 12 in \cite{royden2010real}, page 42), we can find disjoint open intervals $E_k \subset [\frac{1}{4}, z_1], 0\le k\le n$ such that
\begin{equation*}
    \mu(\bigcup_{k=1}^n E_k\setminus E_0)+\mu(E_0\setminus \bigcup_{k=1}^n E_k)<\frac{\epsilon^2}{2}.
\end{equation*}
By the result of Case (ii), we have
\begin{equation*}
    \mathbf{1}_{\bigcup_{k=1}^n E_k} =\sum_{k=1}^n\mathbf{1}_{E_k} \in \mathcal{H},
\end{equation*}
this implies 
\begin{eqnarray*}
  &&\|\mathbf{1}_{\bigcup_{k=1}^n E_k} -\mathbf{1}_{E}\|_2^2 \cr
  &=& \int_{(E_0\setminus\bigcup_{k=1}^n E_k)\bigcup(\bigcup_{k=1}^n E_k\setminus E_0)}(\mathbf{1}_{\bigcup_{k=1}^n E_k} -\mathbf{1}_{E_0})^2\, d\mu \cr
  &&\quad\quad\quad + \int_{[z_1,+\infty)}(\mathbf{1}_{\bigcup_{k=1}^n E_k} -\mathbf{1}_{E})^2\, d\mu\cr
  &=& \mu(\bigcup_{k=1}^n E_k\setminus E_0)+\mu(E_0\setminus \bigcup_{k=1}^n E_k)
  +\mu([z_1,+\infty)) \cr
 &<& \frac{\epsilon^2}{2}+\frac{\epsilon^2}{2} =\epsilon^2,
\end{eqnarray*}
hence $\|\mathbf{1}_{\bigcup_{k=1}^n E_k} -\mathbf{1}_{E_0}\|_2<\epsilon$, which shows $\mathbf{1}_E \in \mathcal{H}$ for $E\in \mathcal{M}$.
 \end{enumerate}
  \end{proof}

\medskip

\emph{The proof of Theorem \ref{thm6}}. It suffices to show $L^2(\mathcal{Z},\mathcal{M},\mu) = \mathcal{H}$. Let $\epsilon$ be arbitrary in $(0,1)$ and $f$ be arbitrary in $L^2(\mathcal{Z},\mathcal{M},\mu)$. 

Firstly, there exists $z_2>\frac{1}{4}$ such that 
\begin{equation}\label{eq54}
    \int_{[z_2,+\infty)}f^2\, d\mu < \frac{\epsilon^2}{4}.
\end{equation}

Secondly, for fixed $f$ and $n\ge 1$, let 
\[
A_{n}^{f}=\left\{z\in \left[\frac{1}{4}, z_2\right]\colon |f(z)|\ge 2^n\right\}, 
\] 
then $A_{n}^{f}\in \mathcal{M}$ and is decreasing as $n\to \infty$ 
. We claim that there exists $N\ge 1$ depending on $\epsilon$ and $z_2$, such that for all $n>N$, there holds 
\begin{equation}\label{eq55}
    \int_{A_{n}^{f}} f^2\, d\mu<\frac{\epsilon^2}{4}.
\end{equation}
We argue it by contradiction. Suppose the claim does not hold. Then there exists some $\epsilon_0>0$ and a subsequence of positive integers $\{n_k\colon k\ge 1\}$ such that \begin{equation*}
    \int_{A_{n_k}^{f}} f^2\, d\mu\ge\frac{\epsilon_0^2}{4}.
\end{equation*}
We denote 
\[
A_{\infty}^{f}= \{z\in \left[\frac{1}{4}, z_2\right]\colon f(z)=\pm \infty\},
\]
then  
\[
\lim_{k\to \infty}A_{n_k}^{f} =\bigcap_{k=1}^{\infty}A_{n_k}^{f} = A_{\infty}^{f}
\]
and 
\begin{equation*}
    \int_{A_{\infty}^{f}} f^2\, d\mu\ge\frac{\epsilon_0^2}{4}.
\end{equation*}
This shows 
\[
\mu(A_{\infty}^{f})>0.
\]
Then we can deduce that
\begin{equation*}
    \int_{\mathcal{Z}} f^2\, d\mu\ge \int_{A_{\infty}^{f}} f^2\, d\mu=+\infty,
\end{equation*}
we have a contradiction since $f\in L^2(\mathcal{Z},\mathcal{M},\mu)$. Therefore the claim holds true.

Thirdly, let 
\[
n_0= \max\left\{N+1,\left\lfloor \log_2\frac{2\sqrt{4z_2-1}}{\sqrt{\pi}\epsilon}\right\rfloor+1\right\}
\]
and  
\begin{equation*}
    f_1(z) = \left\{
\begin{array}{ll}
\displaystyle 2^{-n_0}\lfloor f(z)2^{n_0} \rfloor,& \displaystyle \textrm{if } z \in \left[\frac{1}{4}, z_2\right]\setminus A_{n_0}^{f}, \vspace{2mm}\cr
0,& \textrm{otherwise.}
\end{array}
\right.
\end{equation*}
Then $f_1$ is \emph{a simple function}, i.e. the linear combination of indicators. By employing Lemma \ref{lm11}, we have $f_1\in \mathcal{H}$. In addition, 
\begin{eqnarray}
  \int_{[\frac{1}{4}, z_2]\setminus A_{n_0}^{f}} (f-f_1)^2 \, d\mu&=& \int_{[\frac{1}{4}, z_2]\setminus A_{n_0}^{f}} (f-f_1)^2(z) \frac{\sqrt{4z-1}}{2\pi z^2}\, dz  \nonumber\\
  &<&\int_{[\frac{1}{4}, z_2]\setminus A_{n_0}^{f}} (f-f_1)^2(z) \frac{1}{\pi z^{\frac{3}{2}}}\, dz  \nonumber\\
  &\le&\int_{[\frac{1}{4}, z_2]\setminus A_{n_0}^{f}} 4^{-n_0} \frac{8}{\pi}\, dz  \nonumber\\
  &\le&\int_{[\frac{1}{4}, z_2]} 4^{-n_0} \frac{8}{\pi}\, dz  \nonumber\\
  &\le&\frac{8(z_2-\frac{1}{4})}{\pi}4^{-n_0}   \nonumber\\
  &<&\frac{8(z_2-\frac{1}{4})}{\pi} \frac{\pi \epsilon^2}{16(z_2-\frac{1}{4})}  =\frac{\epsilon^2}{2} \label{eq56}
\end{eqnarray}
Now by employing \eqref{eq54}-\eqref{eq56}, we have
\begin{eqnarray*}
\int_{\mathcal{Z}} (f-f_1)^2 \, d\mu
  &=&\int_{[\frac{1}{4}, z_2]\setminus A_{n_0}^{f}} (f-f_1)^2 \, d\mu
  +\int_{A_{n_0}^{f}} f^2 \, d\mu
  +\int_{[z_2,+\infty)} f^2 \, d\mu \vspace{2mm}\cr
  &<& \frac{\epsilon^2}{2}+\frac{\epsilon^2}{4} +\frac{\epsilon^2}{4} = \epsilon^2,
\end{eqnarray*}
hence 
\[
\|f-f_1\|_2<\epsilon,
\]
 which implies 
 \[
 f-f_1\in \mathcal{H},
 \]
and $f=f_1+(f-f_1)\in \mathcal{H}$. Since $f$ is arbitrary in $L^2(\mathcal{Z},\mathcal{M},\mu)$, we have
 \[L^2(\mathcal{Z},\mathcal{M},\mu) \subset \mathcal{H}.\]
On the other hand, we have also $L^2(\mathcal{Z},\mathcal{M},\mu) \supset \mathcal{H}$ by the definition of $\mathcal{H}$. Thus $L^2(\mathcal{Z},\mathcal{M},\mu) = \mathcal{H}$. \qed

\medskip

Theorem \ref{thm6} means that  $\{g_{2n}\colon n \in Z_+\}$ is an orthonormal basis on semi-infinite interval $\mathcal{Z}$ with respect to the weight function 
\[
w(z) =\frac{2}{\pi}r(z)=\frac{\sqrt{4z-1}}{2\pi z^2},
\] 
where $r(z)$ is defined as in Theorem \ref{thm03}. Other Hilbert spaces and related orthonormal bases  with respect to weight functions, say, four kinds of Chebyshev polynomials, Legendre polynomials, Laguerre polynomials, have shown very powerful influences in orthogonal approximations, scientific computations and vast range of real-world applications(see \cite{jangid2023application},\cite{junghanns2021weighted},\cite{pucanovic2023chebyshev}, \cite{ye2023novel},\cite{yang2023chebyshev},\cite{zhang2023localized}). In Table \ref{table1} we compare $\{g_{2n}\colon n \in Z_+\}$ with some classical orthonormal systems.

\begin{table}[!htbp]
    \centering
    \caption{Some orthonormal systems}
    \begin{tabular}{|c|c|c|c|}
    \hline
    Name &
   $\begin{array}{c}
        \textrm{basis}  \\
         (n\in \Z_+)
   \end{array}$  & domain &$\begin{array}{c}
        \textrm{weight}  \\
         \textrm{function}
   \end{array}$  \\
         \hline
         $\begin{array}{c}
        \textrm{Chebyshev polynomials}  \\
         \textrm{of the first kind}
   \end{array}$  & $\{T_{n}\}$  & $[-1,1]$ &$\frac{1}{\sqrt{1-x^2}}$   \\
         \hline
         $\begin{array}{c}
        \textrm{Chebyshev polynomials}  \\
         \textrm{of the second kind}
   \end{array}$   & $\{U_{n}\}$  & $[-1,1]$ &$\sqrt{1-x^2}$   \\  
         \hline
         normalized Legendre polynomials  & $\{P_{n}\}$ & $[-1,1]$ &1   \\  
         \hline
         Laguerre polynomials  & $\{L_{n}\}$  & $[0,+\infty)$ &$e^{-x}$   \\  
         \hline
         Hermite polynomials  & $\{H_{n}\}$  & $(-\infty, +\infty)$ &$e^{-x^2}$   \\  
         \hline
         $g$-rational functions  & $\{g_{2n}\}$  & $[\frac{1}{4}, +\infty)$ &$\frac{\sqrt{4x-1}}{2\pi x^2}$   \\
         \hline
    \end{tabular}
    \label{table1}
\end{table}
where for $n\in \Z_+$, $g_{2n}$ and $U_n$ are defined as in \eqref{eq43} and \eqref{eq61},  and
\begin{eqnarray*}
  T_n(x)&=&\cos n\arccos x,\\
  P_n(x)&=& \sqrt{\frac{2n+1}{2}}\frac{1}{n!}\sum_{m=0}^{\left\lfloor \frac{n}{2}\right\rfloor}(-1)^m\frac{(2n-2m)!}{m!(n-m)!(n-2m)!}x^{n-2m},\\
  L_n(x) &=& \sum_{m=0}^n (-1)^m \frac{n!}{(n-m)!(m!)^2}x^m,\\
  H_n(x)&=& n!\sum_{m=0}^{\left\lfloor \frac{n}{2}\right\rfloor}(-1)^m\frac{(2x)^{n-2m}}{m!(n-2m)!}.
\end{eqnarray*}

 We see that the sequence of $g$-rational functions differs from other basis from domain and weight function. Both $g$-rational function sequence and Laguerre polynomials sequence are defined on semi-infinite intervals, Hermite polynomials are defined on the entire real line, and the other three kinds are defined on a finite interval $[-1, 1]$. Correspondingly, three weight functions defined on infinite intervals decay when tending to $+\infty$, among which the weight of $g$-rational functions decays ''slowest''. We wish $g$-rational functions could apply to weighted rational approximation examples in future works. For more about rational approximations, please refer to \cite{graves2006rational}.

 A collection $\mathcal{F}$ of normed linear space $X$ is said to be \emph{dense} in $X$ if for each $f\in X$ and $\epsilon>0$, there exists $g\in \mathcal{F}$ such that $\|f-g\|_p<\epsilon$. A normed linear space $X$ is said to be \emph{separable} if there is a countable set dense in $X$. Since the linear combinations of finite number of $g_{2n}$'s with rational coefficients are dense in  $L^2(\mathcal{Z},\mathcal{M},\mu)$, we obtain
\begin{corollary}
$L^2(\mathcal{Z},\mathcal{M},\mu)$ is a separable space.
\end{corollary}

\subsection{A non-orthogonal basis}
 In this part we give a non-orthogonal basis, and prove that it leads to $\{g_{2n}\colon n \in Z_+\}$ after the Gram-Schmidt orthogonalization process. This process is a main approach to obtain an orthonormal basis from a non-orthogonal one. Now we do some preparations on combinatorial numbers.  
 \begin{lemma}\label{lm13}
  Suppose $n\ge m\ge 3$. here holds
\begin{equation}\label{eq59}
    \binom{n}{m} + \binom{n}{m-1} - \binom{n}{m-2} - \binom{n}{m-3} = \binom{n+2}{m} - \binom{n+2}{m-1}.
\end{equation}
 \end{lemma}
\begin{proof}
Recall that for $n\ge m \ge 1$, there holds
\[
\binom{n}{m}+\binom{n}{m-1} =\binom{n+1}{m}, 
\] 
Now let $n\ge m\ge 3$. We have 
\begin{eqnarray*}
    &&\binom{n}{m} + \binom{n}{m-1} - \binom{n}{m-2} - \binom{n}{m-3} \cr
    &=&\binom{n+1}{m}  - \binom{n+1}{m-2} \cr
    &=&\binom{n+1}{m}  + \binom{n+1}{m-1} - \binom{n+1}{m-1}- \binom{n+1}{m-2} \cr
    &=&\binom{n+2}{m}  - \binom{n+2}{m-1}.
\end{eqnarray*}
\end{proof}
In the following proposition, we represent $z^{-n}$ in terms of $\{g_{2n}\colon n\in \Z_+\}$.
\begin{proposition}\label{prop6}
  For all $n\in\mathbb{Z}_+$, $z\in \C^*$, the following identities hold true:
\begin{equation}\label{eq09}
  z^{-n}= g_{2n}(z)+\sum_{l=0}^{n-1} \left[\binom{2n}{n-l}-\binom{2n}{n-l-1}\right]g_{2l}(z).
\end{equation}
\end{proposition}
\begin{proof}
Since for each $z\in \C^*$,
\begin{equation*}
    \left\{
\begin{array}{lcl}
   g_0(z)  &\equiv& 1, \\
    \displaystyle g_2(z) & =& z^{-1} - 1,\\
    g_4(z) &=& z^{-2}-3z^{-1}+1,\\
     g_6(z)&= & z^{-3} - 5z^{-2}+ 6z^{-1} - 1,
\end{array}
    \right.
\end{equation*}
Solving this linear system, we obtain
\begin{equation*}
\left\{
\begin{array}{lcl}
   z^0  &=&  g_0(z), \\
    z^{-1} & =& g_0(z)+g_2(z) \\
    z^{-2} &=&  g_4(z) +3g_2(z) +2g_0\\
     z^{-3}&= & g_6(z)+5g_4(z) +9g_2(z) +5g_0.
\end{array}
\right.
\end{equation*}
That is,
\begin{equation*}
\left\{
\begin{array}{lcl}
   z^0  &=& g_0(z), \vspace{2mm}\\
    z^{-1} & =& g_2(z)+\left[\binom{2}{1}-\binom{2}{0}\right]g_0(z),
    \vspace{2mm}\\
    z^{-2} &=&  \vspace{2mm}g_4(z)+\left[\binom{4}{1}-\binom{4}{0}\right]g_2(z)+\left[\binom{4}{2}-\binom{4}{1}\right]g_0(z)\\
     z^{-3}&= &  g_6(z)+\left[\binom{6}{1}-\binom{6}{0}\right]g_4(z)+\left[\binom{6}{2}-\binom{6}{1}\right]g_2(z)+\left[\binom{6}{3}-\binom{6}{2}\right]g_0(z).
\end{array}
\right.
\end{equation*}
Hence \eqref{eq09} holds for $n\le 3$. Now let $n\ge4$ and suppose  \eqref{eq09} holds for $n-1$. That is, for $z\in \C^*$,
\begin{eqnarray}\label{eq60}
  z^{-n+1}&= &g_{2n-2}(z)+\sum_{l=0}^{n-2} \left[\binom{2n-2}{n-l-1}-\binom{2n-2}{n-l-2}\right]g_{2l}(z).
\end{eqnarray}

By Lemma \ref{lm02}, Identity \eqref{eq69}, we have for $n\in \Z_+, z\in \C^*$, 
\begin{equation*}
  g_{2n+4}(z) = \left(\frac{1}{z}-2\right)g_{2n+2}(z)-g_{2n}(z),
\end{equation*}
This implies
\begin{equation}\label{eq62}
  \frac{1}{z}g_{2n+2}(z)= g_{2n+4}(z) + 2g_{2n+2}(z)+g_{2n}(z).
\end{equation}
Using equation \eqref{eq60}, we may write 
\begin{equation}\label{eq63}
  z^{-n} = \frac{1}{z} g_{2n-2}(z)+\sum_{l=0}^{n-2} \left[\binom{2n-2}{n-l-1}-\binom{2n-2}{n-l-2}\right]\frac{1}{z}g_{2l}(z).
\end{equation}
By employing \eqref{eq62} and \eqref{eq63}, we have

\begin{eqnarray*}
  &&z^{-n}\\
  &=& z^{-1}z^{-n+1} \\
  &=& [g_{2n}(z) + 2g_{2n-2}(z)+g_{2n-4}(z)]+\sum_{l=1}^{n-2} \left[\binom{2n-2}{n-l-1}-\binom{2n-2}{n-l-2}\right]\\
  && [g_{2l+2}(z) + 2g_{2l}(z)+g_{2l-2}(z)] + \left[\binom{2n-2}{n-1}-\binom{2n-2}{n-2}\right][g_0(z)+g_2(z)]\\
  &=& [g_{2n}(z) + 2g_{2n-2}(z)+g_{2n-4}(z)]+\sum_{l=1}^{n-2} \left[\binom{2n-2}{n-l-1}-\binom{2n-2}{n-l-2}\right]\\
   &&g_{2l+2}(z) +\sum_{l=1}^{n-2} 2\left[\binom{2n-2}{n-l-1}-\binom{2n-2}{n-l-2}\right]
   g_{2l}(z) \\
   &&+\sum_{l=1}^{n-2} \left[\binom{2n-2}{n-l-1}-\binom{2n-2}{n-l-2}\right]
   g_{2l-2}(z)   \\
  && + \left[\binom{2n-2}{n-1}-\binom{2n-2}{n-2}\right][g_0(z)+g_2(z)].
  \end{eqnarray*}
  Changing the index, we obtain
\begin{eqnarray*}
  z^{-n} &=& [g_{2n}(z) + 2g_{2n-2}(z)+g_{2n-4}(z)]\\
  &&+\sum_{l'=2}^{n-1} \left[\binom{2n-2}{n-l'}-\binom{2n-2}{n-l'-1}\right]
   g_{2l'}(z) \\
   &&+\sum_{l=1}^{n-2} \left[2\binom{2n-2}{n-l-1}-2\binom{2n-2}{n-l-2}\right]
   g_{2l}(z) \\
   &&+\sum_{l''=0}^{n-3} \left[\binom{2n-2}{n-l''-2}-\binom{2n-2}{n-l''-3}\right]
   g_{2l''}(z)\\
   &&+ \left[\binom{2n-2}{n-1}-\binom{2n-2}{n-2}\right][g_0(z)+g_2(z)]
   \end{eqnarray*}
  Separating $g_{2n}, g_{2n-2}, g_{2n-4}, g_{2}$ and $g_0$ terms and rearrange the rest terms, we get
\begin{eqnarray}
 z^{-n} &=& [g_{2n}(z) + (2n-1)g_{2n-2}(z)+(2n^2-3n)g_{2n-4}(z)]\nonumber\\
 &&+\sum_{l=2}^{n-3} \left[\binom{2n-2}{n-l}+\binom{2n-2}{n-l-1}-\binom{2n-2}{n-l-2}-\binom{2n-2}{n-l-3}\right]g_{2l}(z) \nonumber\\
   &&+ \left[ \binom{2n-2}{n-1}+\binom{2n-2}{n-2} - \binom{2n-2}{n-3} -\binom{2n-2}{n-4} \right]g_2(z)\nonumber\\
   &&+ \left[\binom{2n-2}{n-1}-\binom{2n-2}{n-3}\right]g_0(z).\label{eq64}
\end{eqnarray}
Since $\displaystyle \binom{2n-1}{n}=\binom{2n-1}{n-1}$, we have  
\begin{eqnarray}
  \binom{2n-2}{n-1}-\binom{2n-2}{n-3} &=& \binom{2n-2}{n-1}+\binom{2n-2}{n-2}-\binom{2n-2}{n-2}-\binom{2n-2}{n-3}\nonumber\\
  &=& \binom{2n-1}{n-1}-\binom{2n-1}{n-2}\nonumber\\
  &=& \binom{2n-1}{n-1}+\binom{2n-1}{n}-\binom{2n-1}{n-1}-\binom{2n-1}{n-2}\nonumber\\
  &=& \binom{2n}{n}-\binom{2n}{n-1}\label{eq65}.
\end{eqnarray}
Applying \eqref{eq59} and \eqref{eq65} in the right hand side of \eqref{eq64}, we obtain
\begin{equation*}
  z^{-n} =g_{2n}(z)+\sum_{l=0}^{n-1} \left[\binom{2n}{n-l}-\binom{2n}{n-l-1}\right]g_{2l}(z).
\end{equation*}
By induction, equation \eqref{eq09} holds for all $n\in \Z_+$.
\end{proof}
We have the following theorem
\begin{theorem}\label{thm9}
\begin{enumerate}[(i)]
    \item $\{z^{-n}\colon n \in \Z_+\}$ forms a basis for $L^2(\mathcal{Z},\mathcal{M},\mu)$.
    \item The Gram-Schmidt process of $\{z^{-n}\colon n \in \Z_+\}$ in $L^2(\mathcal{Z},\mathcal{M},\mu)$  obtains $\{g_{2n}\colon n \in Z_+\}$.
    \item For each $n\in \Z_+$, $\mathrm{span}\{z^{-k}\colon 0\le k\le n\} = \mathrm{span}\{g_{2k}\colon 0\le k\le n\}$.
\end{enumerate}
\end{theorem}
\begin{proof}
\begin{enumerate}[(i)]
    \item For all $n\in \Z_+$,  Proposition \ref{prop6} implies $z^{-n}\in L^2(\mathcal{Z},\mathcal{M},\mu)$. Now suppose $f \in L^2(\mathcal{Z},\mathcal{M},\mu)$. Given $\epsilon>0$, by Theorem \ref{thm6} there exists some $n$ and $\{c_k\colon 0\le k\le n\}$ such that 
\[
\|f-\sum_{k=1}^{n}c_mg_{2m}\|_2 = \|f-\sum_{m=1}^{n}c_k\sum_{l=0}^{m}(-1)^l\binom{m-l}{l}z^{-(m-l)}\|_2< \epsilon.
\]
This shows 
$$f\in \overline{\mathrm{span}}\{z^{-n}\colon n\in \Z_+\}.$$ Hence  $\{z^{-n}\colon n \in \Z_+\}$ is a basis for $L^2(\mathcal{Z},\mathcal{M},\mu)$. 
\item We use the symbol $v_{n}$ to denote the $(n+1)$-th orthogonal function obtained by applying the Gram-Schmidt process to $\{z^{-j}\colon \in \Z_+\}$. Then $v_0 = z^0 =g_0$. For each $n\in \N$, suppose that $v_k =g_{2k}, k<n$. Then by \eqref{eq09}, we can derive
\begin{eqnarray*}
  v_n&=&z^{-n}-\sum_{k=0}^{n-1}\frac{(z^{-n},v_k) v_k}{(v_k,v_k)}\\
  &=&z^{-n}-\sum_{k=0}^{n-1}(z^{-n},g_{2k}) g_{2k}\\
  &=& z^{-n}-\sum_{k=0}^{n-1}\left( g_{2n}(z)+\sum_{l=0}^{n-1} \left[\binom{2n}{n-l}-\binom{2n}{n-l-1}\right]g_{2l}(z) ,g_{2k}\right)g_{2k}\\
  &=& z^{-n}-\sum_{k=0}^{n-1}\left[\binom{2n}{n-l}-\binom{2n}{n-l-1}\right]g_{2k}
  =g_{2n}.
\end{eqnarray*}
\item It is clear by applying \eqref{eq43} and \eqref{eq09}.
\end{enumerate}
\end{proof}
The following corollary is about the transition matrices of the two bases. 
\begin{corollary}\label{cor06}
For all $n\in \N$, the following identity holds true:
\[A_nB_n= I_{n},\] 
where $I_{n}$ is the order $n$ identify matrix, $A_n$ and $B_n$ are both order $n$ lower triangular matrices, $\displaystyle A_n=(a_{i,j})_{1\le i,j\le n}, 
B_n=(b_{i,j})_{1\le i,j\le n}$, and for $1\le i,j\le n$,
\begin{eqnarray*}
  a_{ij}&=&\left\{
\begin{array}{ll}
\displaystyle (-1)^{i-j}\binom{i+j-2}{i-j},     & \textrm{if } i\ge j; \\
   0,  & \textrm{otherwise}
\end{array}
\right.,\\
b_{ij}&=&\left\{
\begin{array}{ll}
\displaystyle\binom{2i-2}{i-j}- \binom{2i-2}{i-j-1},    & \textrm{if } i> j; \\
1,    & \textrm{if } i= j; \\
   0,  & \textrm{otherwise}.
\end{array}
\right.
\end{eqnarray*}
\end{corollary}
\begin{proof}
$A_n$ is the transition matrix from $\{z^{-k}\colon 0\le k\le n-1\}$ to $\{g_{2k}\colon 0\le k\le n-1\}$ and $B_n$ is the one from $\{g_{2k}\colon 0\le k\le n-1\}$ to $\{z^{-k}\colon 0\le k\le n-1\}$, hence $B_n$ is the inverse of $A_n$.
\end{proof}

 For example, when $n=5$, 
\[
A_5=\begin{bmatrix}
1&&&&\\
-1&1&&&\\
1&-3&1&&\\
-1&6&-5&1&\\
1&-10&15&-7&1
\end{bmatrix}, \quad B_5=\begin{bmatrix}
1&&&&\\
1&1&&&\\
2&3&1&&\\
5&9&5&1&\\
14&28&20&7&1
\end{bmatrix},
\]
Thus we have $A_5B_5=I_5$.

\subsection{Examples} 
  For convenience, In the sequel we will call a Fourier expansion of function $f$ in with respect to $\{g_{2n}\colon n \in Z_+\}$ in $L^2(\mathcal{Z},\mathcal{M},\mu)$ \emph{a Fourier expansion of function $f$} for short. In this part we give two examples as applications, including one for Fourier series and one for interpolation. First we expand function  $g_{2m+1}(z) = \frac{1}{\sqrt{z}}$ on  $\mathcal{Z}$ using basis $\{g_{2n}\colon  n\in \Z_+\}$ for all $m\in \Z_+$ in the following example. 
 \begin{example}
 For all $m\in \Z_+$, 
 \begin{enumerate}[(i)]
     \item $\|g_{2m+1}\|=1$.
     \item The Fourier series expansion of $g_{2m+1}$ is 
 \begin{equation}\label{eq58}
     g_{2m+1}=  \sum_{n=0}^{\infty} \frac{8(m+1)(-1)^{n+m+1}g_{2n}}{\pi(2n-2m-1)(2n+2m+3)}.
 \end{equation}
 \end{enumerate}
 \end{example}
  \begin{proof}
  For each $m\in \Z_+$, by applying  \eqref{eq72} and \eqref{eq51}, we have
  \begin{equation*}
    \|g_{2m+1}\|^2 =\int_{\mathcal{Z}}g_{2m+1}^2 \, d\mu =\frac{4}{\pi}\int_{0}^{\frac{\pi}{2}}\sin^2(2m+2)\theta \, d\theta =1.
  \end{equation*}
  Hence $\|g_{2m+1}\|=1$. For each $m,n\in \Z_+$, by applying Lemma \ref{lm12} and \eqref{eq51}, we deduce that 
  \begin{eqnarray*}
   (g_{2m+1}, g_{2n}) &=& \int_{\mathcal{Z}} g_{2m+1} g_{2n}\, d\mu \\
   &=& \frac{4}{\pi}\int_{0}^{\frac{\pi}{2}} \frac{\sin(2m+2)\theta}{\sin \theta} \frac{\sin (2n+1)\theta}{\sin \theta}  \sin^2\theta\, d\theta \\
   &=& \frac{2}{\pi} \left[\frac{\sin(2n-2m-1)\theta}{2n-2m-1}  - \frac{\sin(2n+2m+3)\theta}{2n+2m+3}\right]_0^{\frac{\pi}{2}} \\
   &=& \frac{8(m+1)(-1)^{n+m+1}}{\pi(2n-2m-1)(2n+2m+3)}.
 \end{eqnarray*}
Hence in the closed form, The Fourier series of $g_{2m+1}$ is
 \begin{equation*}
     g_{2m+1}=  \sum_{n=0}^{\infty} \frac{8(m+1)(-1)^{n+m+1}g_{2n}}{\pi(2n-2m-1)(2n+2m+3)}. \end{equation*}
 \end{proof}
 We have the following identity.
 \begin{proposition}
Let $\mathbb{Z}^* = \Z\setminus\{0\}$.  For each $m\in \mathbb{Z}^*$, the following identity holds true:
   \begin{equation}\label{eq57}
     \frac{\pi^2}{64}=  \sum_{n=0}^{\infty} \frac{m^2}{(2n-2m+1)^2(2n+2m+1)^2}.
 \end{equation}
 \end{proposition}
 \begin{proof}
 For $m\in \Z_+$, since $\|g_{2m+1}\|=1$, by the Parseval's identity of $g_{2m+1}$'s and equation \eqref{eq58}, we have
 \begin{equation*}
     1=  \sum_{n=0}^{\infty} \frac{64(m+1)^2}{\pi^2(2n-2m-1)^2(2n+2m+3)^2}.
 \end{equation*}
 Multiply both sides of equation above by $\frac{\pi^2}{64}$ and let $m'= m+1$,  we have for $m\ge 1$, \eqref{eq57} holds true. Meanwhile, it is true as the same if we replace $m$ by $- m$. Thus, it holds true for all $m\in \mathbb{Z}^* $.
 \end{proof}
For example, the Fourier expansion of $g_{1}$  is as follows. 
 \begin{eqnarray*}
     g_1 & = &  \sum_{n=0}^{\infty} \frac{8(-1)^{n+1}g_{2n}}{\pi(2n-1)(2n+3)}\\
     &=&\frac{8}{\pi}\left[-\frac{1}{(-1)\cdot3}g_0+\frac{1}{1\cdot5}g_2-\frac{1}{3\cdot 7}g_4+\frac{1}{5\cdot 9}g_6+\cdots\right].
 \end{eqnarray*}
 The order $n$ best approximation of $\frac{1}{\sqrt{z}}$ in the sense of $L^2(\mathcal{Z},\mathcal{M},\mu)$, denoted by $P_n\left(\frac{1}{\sqrt{z}}\right)$ is
 \[
P_n\left(\frac{1}{\sqrt{z}}\right) = \frac{8}{\pi}\left[\frac{1}{3}g_0+\frac{1}{1\cdot5}g_2-\frac{1}{3\cdot 7}g_4+\frac{1}{5\cdot 9}g_6+\cdots +  \frac{(-1)^{n+1}g_{2n}}{(2n-1)(2n+3)}\right].
 \]
 
 \begin{figure}[tbhp]
	\centering  
	\subfigure{
	\includegraphics[width=0.48\linewidth]{"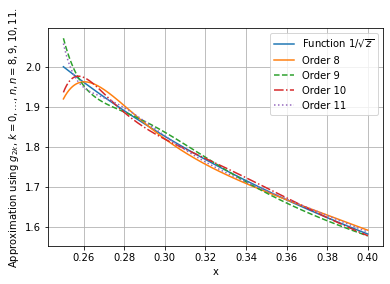"}}
	\subfigure{
	\includegraphics[width=0.48\linewidth]{"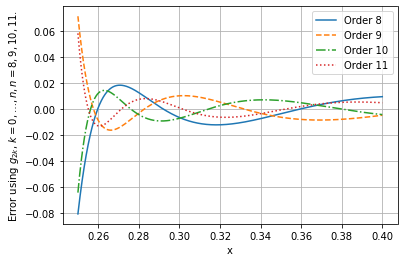"}}
	\caption{Approximation and error curves of partial sums of $g_{0}. g_{2}, \dots, g_{22}$ in approximating $\frac{1}{\sqrt{z}}$}
	\label{fig: interpo}
\end{figure}
The best approximation curves and corresponding error curves for cases  $n=8$ to $n=11$ on the interval $[\frac{1}{4}, \frac{2}{5}]$ are shown in Figure \ref{fig: interpo}.
 The left sub-figure shows the approximation curves, the right sub-figure shows the error curves. The corresponding $L^2(\mathcal{Z},\mathcal{M},\mu)$ errors are 
\begin{eqnarray*}
  e_8 &=&0.013671986780782784, \cr
e_9 &=&0.011663773398010387,\cr
e_{10} &=&0.010103838132371137,\cr
e_{11} &=&0.00886345373508862.
\end{eqnarray*}

In the following we give a simple interpolation example.
 \begin{example}
 Suppose 
 \[
 f(1) =2, f(2) =1, f(3)=1.5, f(4) =1.
 \] Now we determine a function in $\overline{\textrm{span}}\{g_0,g_2,g_4,g_6\}$ that interpolates these four points. Let
the interpolant 
\[
Pf = c_1g_0+c_2g_2+c_3g_4+c_4g_6,
\]
then in matrix form we have 
\begin{equation*}
    \begin{bmatrix}
    g_0(1)& g_2(1)& g_4(1)& g_6(1)\\
    g_0(2)& g_2(2)& g_4(2)& g_6(2)\\
    g_0(3)& g_2(3)& g_4(3)& g_6(3)\\
    g_0(4)& g_2(4)& g_4(4)& g_6(4)
    \end{bmatrix}
    \begin{bmatrix}
    c_1\\
    c_2\\
    c_3\\
    c_4
    \end{bmatrix}
    = \begin{bmatrix}
    y_1\\
    y_2\\
    y_3\\
    y_4
    \end{bmatrix},
\end{equation*}
Substituting the values $g_{2n}(j)$ into the system above, $n=0,\dots,3, j=1,\dots, 4$, we obtain
\begin{equation*}
    \begin{bmatrix}
   1& 0& -1&1\vspace{2mm}\\
    1& \frac{1}{2}& -\frac{1}{4}& \frac{7}{8}\vspace{2mm}\\
    1& \frac{2}{3}& \frac{1}{9}& \frac{13}{27}\vspace{2mm}\\
    1& \frac{3}{4}& \frac{5}{16}& \frac{13}{64}
    \end{bmatrix}
    \begin{bmatrix}
    c_1\vspace{2mm}\\
    \vspace{2mm}c_2\\
    \vspace{2mm}c_3\\
    c_4
    \end{bmatrix}
    = \begin{bmatrix}
    2\vspace{2mm}\\
    1\vspace{2mm}\\
    \vspace{2mm}\frac{3}{2}\\
    1
    \end{bmatrix},
\end{equation*}
By solving this linear system, we get
\begin{equation*}
    \begin{bmatrix}
    \vspace{2mm}c_1\\
    \vspace{2mm}c_2\\
    \vspace{2mm}c_3\\
    c_4
    \end{bmatrix}
    = \begin{bmatrix}
    \frac{811}{6}\vspace{2mm}\\
    \frac{1097}{4}\vspace{2mm}\\
    \frac{1147}{6}\vspace{2mm}\\
    \frac{70}{113}
    \end{bmatrix} \approx \begin{bmatrix}
    \vspace{2mm}14.072271386430678\\
   \vspace{2mm} -18.825221238938052\\
    \vspace{2mm}11.157817109144542\\
    \vspace{2mm}-0.6194690265486725
    \end{bmatrix}.
\end{equation*}
Hence, the interpolant $Pf$ is
\[
Pf(z) = \frac{811}{6}g_0(z)+ \frac{1097}{4}g_2{z} +\frac{1147}{6}g_4(z)+ \frac{70}{113}g_6(z),\quad z\in \C^*.
\]
 \end{example}

\section{Conclusion}
In this article , Through investigating four sequences, including a polynomial sequence $\{G_n\colon n\in \Z_+\}$, an irrational function sequence $\{g_n\colon n\in \Z_+\}$, two rational function sequences $\{g_{2n}\colon n\in \Z_+\}$ and $\{z^{-n}\colon n\in \Z_+\}$, we established a Hilbert space $L^2(\mathcal{Z},\mathcal{M},\mu)$. We are intended to reveal more properties of this space in the future work.

\section*{Acknowledgments}
The author extends his gratitude to Xingping Sun from Missouri State University for his insightful guidance and earnest help. 
The author would like to thank the anonymous referees for their valuable suggestions.

\begin{thebibliography}{99}
\bibitem{frontczak2023balancing}
R. Frontczak and K. Prasad. Balancing polynomials, fibonacci numbers and some series for $\pi$, {\it Mediterranean Journal of Mathematics, } {\bf 20}(4), 1--18(2023).

\bibitem{marques2022proof}
D. Marques and P. Trojovsk{\`y}. The proof of a formula concerning the asymptotic behavior of the reciprocal sum of the square of multiple-angle fibonacci numbers, {\it Journal of inequalities and Applications, } {\bf 2022}(1), 21(2022).

\bibitem{suvarnamani2015some}
A. Suvarnamani and M. Tatong. Some properties of $(p,q)$-fibonacci numbers,  {\it Progress in Applied Science and Technology, } {\bf 5}(2), 17--21(2015).

\bibitem{kilmer2023constructing}
S. Kilmer, J. Liu, X. Sun and M. Wright. Constructing positive definite polynomials,  {\it Proceedings of the American Mathematical Society, } {\bf 151}(10), 4307--4316(2023).

\bibitem{brualdi2004introductory}
  
  R. A. Brualdi. {\it Introductory Combinatorics, 4th edition,} (2004).
  
  \bibitem{szego1975orthogonal} 
  G. Szeg{\"o}. {\it Orthogonal Polynomials}, Vol. 23, in: American Mathematical Society Colloquium Publications, 1975.
  
  \bibitem{asmar2016partial}
  
  N. H. Asmar. {\it Partial differential equations with Fourier series and boundary value problems,} Courier Dover Publications (2016).
  
  \bibitem{kiepiela2005extension}
  
  K. Kiepiela and D. Klimek.  An extension of the chebyshev polynomials, {\it Journal of Computational and Applied Mathematics,} {\bf 178}(1-2), 305--312 (2005).
  
  \bibitem{BOYD198763}
  J. P. Boyd. Orthogonal rational functions on a semi-infinite interval,  {\it Journal of Computational Physics, } {\bf 70}(1), 63--88(1987).
  
\bibitem{rudin1976principles}
W. Rudin. Principles of mathematical analysis, 3rd edition, (1976).

\bibitem{royden2010real}
  
  H. L. Royden and P. M. Fitzpatrick. {\it Real Analysis, 4th edition,}  Printice-Hall Inc, Boston (2010).

\bibitem{conway2019course}
  
 J. B. Conway. {\it A course in functional analysis}, Vol 96, Springer, (2019).

 \bibitem{stein2011functional}
  
 E. M. Stein and R. Shakarchi. {\it Functional analysis: introduction to further topics in analysis}, Vol 4, Princeton University Press, (2011).

 \bibitem{mason2002chebyshev}
  
 J. C. Mason and D.C. Handscomb. {\it Chebyshev Polynomials.}  CRC Press, (2002).
 
\bibitem{jangid2023application}

K. Jangrid and S. Mukhopadhyay. Applications of legendre wavelet collocation method to the analysis of poro-thermoelastic coupling with variable thermal conductivity, {\it Computers \& Mathematics with Applications,} {\bf 146}, 1--11 (2023).

\bibitem{junghanns2021weighted}
  
P. Junghanns, G. Mastroianni and I. Notarangelo. {\it Weighted polynomial approximation and numerical methods for integral equations}, (2021).

\bibitem{pucanovic2023chebyshev}

Z. Pucanovi{\'c} and M. Pe{\v{s}}ovi{\'c}. Chebyshev polynomials and $r$-circulant matrices, {\it Applied Mathematics and Computation,} {\bf 437}, 127521, (2023).

\bibitem{ye2023novel}

J. Ye, L. Xie, L. Ma, Y. Bian and X. Xu. A novel hybrid model based on Laguerre polynomial and multi-objective Runge--Kutta algorithm for wind power forecasting, {\it International Journal of Electrical Power \& Energy Systems,} {\bf 146}, 108726, (2023).

\bibitem{yang2023chebyshev}

J. Yang, W. Yu, W. Chen, B. Liao and H. Zhu. Chebyshev-Series Solutions for Nonlinear Systems with Hypersonic Gliding Trajectory Example, {\it Aerospace Science and Technology,} 108424, (2023).

\bibitem{zhang2023localized}

S. Zhang and Y. Nie. Localized Chebyshev and MLS collocation methods for solving 2D steady state nonlocal diffusion and peridynamic equations, {\it Mathematics and Computers in Simulation,} {\bf 206}, 264--285, (2023).

\bibitem{graves2006rational}
  
P. R. Graves-Morris, E. B. Saff and R. S. Varga. {\it Rational Approximation and Interpolation: Proceedings of the United Kingdom-United States Conference, Held at Tampa, Florida, December 12-16, 1983}, Vol 1105, (2006).

  \end{thebibliography}

\end{document}